\documentclass[journal]{IEEEtran}
\usepackage{cite}
\usepackage{array}
\usepackage{amsmath,amssymb,dsfont}
\usepackage{amsthm}
\usepackage{verbatim}
\usepackage{cite}
\usepackage{xcolor}
\usepackage{subfigure}
\usepackage{hyperref}
\usepackage{booktabs}
\usepackage[ruled,vlined]{algorithm2e}
\usepackage{graphicx,pbox}
\graphicspath{{figures/}}

\newif\ifarx
 \arxtrue 

\newtheorem{theorem}{Theorem}
\newtheorem{assumption}{Assumption}
\newtheorem{definition}{Definition}
\newtheorem{lemma}{Lemma}

\theoremstyle{definition}
\newtheorem{example}{Example}
\allowdisplaybreaks[1]

\DeclareMathOperator{\dom}{dom}

\DeclareMathOperator*{\argmin}{\arg\min}

\DeclareMathOperator{\T}{\mathsf{T}}
\DeclareMathOperator{\E}{\mathds{E}}
\DeclareMathOperator{\s}{\boldsymbol{s}}
\DeclareMathOperator{\w}{\boldsymbol{w}}

\DeclareMathOperator{\cw}{{\scriptstyle\mathcal{W}}}
\DeclareMathOperator{\bcw}{{\boldsymbol{\scriptstyle\mathcal{W}}}}
\DeclareMathOperator{\cx}{{{\scriptstyle\mathcal{X}}}}
\DeclareMathOperator{\cy}{{{\scriptstyle\mathcal{Y}}}}

\begin{document}
\title{Regularized Diffusion Adaptation via Conjugate Smoothing}

\author{Stefan~Vlaski,~\IEEEmembership{Student Member,~IEEE,} Lieven Vandenberghe,
  and~Ali~H.~Sayed,~\IEEEmembership{Fellow,~IEEE}%
  \thanks{S. Vlaski and A. H. Sayed are with the Institute of Electrical Engineering, \'{E}cole Polytechnique F\'{e}d\'{e}rale de Lausanne. S. Vlaski and L. Vandenberghe are with the Department of Electrical Engineering, University of California, Los Angeles. Emails:\{stefan.vlaski, ali.sayed\}@epfl.ch, vandenbe@ucla.edu. This work was supported in part by NSF grant CIF-1524250. A preliminary version of this work appears in~\cite{Vlaski16}.}
}%

\maketitle

\begin{abstract}
The purpose of this work is to develop and study a distributed strategy for Pareto optimization of an aggregate cost
consisting of regularized risks. Each risk is modeled as the expectation of some loss function with unknown probability
distribution while the regularizers are assumed deterministic, but are not required to be differentiable or even
continuous. The individual, regularized, cost functions are distributed across a strongly-connected network of agents
and the Pareto optimal solution is sought by appealing to a multi-agent diffusion strategy. To this end, the
regularizers are smoothed by means of infimal convolution and it is shown that the Pareto solution of the approximate,
smooth problem can be made arbitrarily close to the solution of the original, non-smooth problem. Performance bounds
are established under conditions that are weaker than assumed before in the literature, and hence applicable to a
broader class of adaptation and learning problems.
\end{abstract}
\IEEEpeerreviewmaketitle%
\section{Introduction}\label{sec:intro}\noindent
The objective of distributed learning is the solution of global, stochastic optimization problems across
networks of agents through localized interactions and without information about the statistical properties of the data.
Using streaming data, the resulting strategies are adaptive in nature and able to track drifts in the location of the minimizers due to variations in the statistical properties of the data. Regularization is one
useful technique to encourage or enforce structural properties on the sought after minimizer, such as sparsity or
constraints. A substantial number of regularizers are inherently non-smooth, while many cost functions are
differentiable. These article proposes a \emph{fully-decentralized} and \emph{adaptive} strategy that is able to
minimize an aggregate sum of regularized costs. To do so, we fully exploit the structure of the individual objectives
as sums of differentiable costs and non-differentiable regularizers.

\emph{Notation:} Throughout the manuscript, random quantities are denoted in boldface. Matrices are denoted in capital letters while vectors and scalars are denoted
in small-case letters. The symbol \( \le \) denotes a regular inequality, while \( \preceq \) denotes an \emph{element-wise} inequality. {Unless specified otherwise, \( \|\cdot\| \) denotes the Euclidean norm.}

\subsection{Problem Formulation}\label{sec:prob_form}\noindent
We consider a strongly-connected network consisting of \(N\) agents. For any two agents \(k\) and \( \ell \), we attach
a pair of  non-negative coefficients \( \{a_{\ell k}, a_{k\ell}\} \) to the edge linking them. The scalar \(a_{\ell k}\)
is used to scale data moving from agent \( \ell \) to \(k\); likewise, for \(a_{k\ell}\). Strong-connectivity means that
it is always possible to find {a path in each direction} with nonzero scaling weights linking any two agents (either
directly if they are neighbors or indirectly through other agents). In addition, at least one agent \(k\) in the network
possesses a self-loop with \(a_{kk}>0\). This condition ensures that at least one agent in the network has some
confidence in its local information. Let \(\mathcal{N}_k\) denote the set of neighbors of agent \(k\). The coefficients
\( \{a_{\ell k}\} \) are convex combination weights that satisfy
\begin{equation}\label{eq:combinationcoef}
  a_{\ell k} \geq 0, \quad \sum_{\ell \in \mathcal{N}_k} a_{\ell k}=1, \quad a_{\ell k} = 0\ \mathrm{if}\ \ell \notin \mathcal{N}_k.
\end{equation}
If we introduce the combination matrix \(A=[a_{\ell k}]\), it then follows from~\eqref{eq:combinationcoef} and the
strong-connectivity property that \(A\) is a left-stochastic primitive matrix. In view of the Perron-Frobenius
Theorem~\cite{Horn03,Pillai05,Sayed14}, this ensures that \(A\) has a single eigenvalue at one while all other
eigenvalues are inside the unit circle, so its spectral radius is given by \( \rho(A)=1 \).  Moreover, if we let \(p\) denote the right-eigenvector
of \(A\) that is associated with the eigenvalue at one, and if we normalize the entries of \(p\) to add up to one, then
it also holds that all entries of \(p\) are strictly positive, i.e.,
\begin{equation}
  Ap=p, \quad \mathds{1}^{\T} p=1, \quad p_k > 0
\end{equation}
where the \( \{ p_k \} \) denote the individual entries of the Perron vector, \(p\).

We associate with each agent \(k\) a risk function \(J_k(w):\mathds{R}^{M} \rightarrow \mathds{R}\), assumed
differentiable. In most adaptation and learning problems, risk functions are expressed as the expectation of loss
functions. Hence, we assume that each risk function is of the form \(J_k(w)=\E Q(w;\boldsymbol{x})\), where \(Q(\cdot)\)
is the loss function and \( \boldsymbol{x} \) denotes random data. The expectation is computed over the distribution of
this data (note that, in our notation, we use boldface letters for random quantities and normal letters for
deterministic quantities or data realizations). We also associate with agent \(k\) a regularization term,
\(R_k(w):\mathds{R}^{M}\rightarrow\mathds{R}\), which is a known deterministic function although possibly
non-differentiable. Regularization factors of this form can, for example, help induce sparsity properties (such as using
\(\ell_1\) or elastic-net regularizers)~\cite{Friedman09,Donoho06,Tibshirani94}.

The objective we are interested in is to devise a fully distributed strategy to seek the {unique} minimizer of the following {strongly-convex}, weighted, aggregate cost, denoted by \(w^o\):
\begin{equation}\label{eq:wnot}
  w^o = \argmin_{w \in \mathds{R}^M} \: \sum_{k=1}^{N} p_k \left \{ J_k(w) + R_k(w)\right \}
\end{equation}
The weights \( \{ p_k \} \) indicate that the resulting minimizer \(w^o\) can be interpreted as a Pareto solution for
the collection of regularized risks \( \{ J_k(w)+R_k(w) \} \)~\cite{Sayed14,Chen13} and will depend on the entries of the Perron eigenvector in a manner specified further below. We are particularly interested in
determining this Pareto solution in the \emph{stochastic} setting when the distribution of the data \(\boldsymbol{x}\)
is unknown. This means that the risks \(J_k(w)\), or their gradient vectors, are also unknown. As such,
\emph{approximate} gradient vectors will need to be employed. A common construction in stochastic approximation theory
is to employ the following choice at each iteration \(i\):
\begin{equation}\label{eq:stochastic_gradient}
  \widehat{\nabla  J}_k(w)\:=\: \nabla  Q_k(w;\boldsymbol{x}_i)
\end{equation}
where \(\boldsymbol{x}_i\) represents the data that is available (observed) at time \(i\). The difference between the
true gradient vector and its approximation is called gradient noise. This noise will seep into the operation of the
distributed algorithm and one main challenge is to show that, despite its presence, the proposed solution is able to
approach \(w^o\) asymptotically. A second challenge we face in constructing an effective distributed solution is the
non-smoothness (non-differentiability) of the regularizers. Motivated by a technique proposed in~\cite{Nesterov05} in
the context of \emph{single} agent optimization, we will address this difficulty in the multi-agent case by introducing
a smoothed version of the regularizers and then showing that the solution \(w^o\) can still be recovered under this
substitution as the size of the smoothing parameter is reduced. We adopt a general formulation that will be shown to
include proximal iterations as a special case.

\subsection{Related Works in the Literature}
\noindent The literature on distributed optimization is extensive. Some early strategies include incremental~\cite{Bertsekas97incremental}, consensus or decentralized gradient descent~\cite{Nedic09, Ram10, Nedic10, Yuan16}, and the diffusion algorithm~\cite{Lopes08, Sayed14, Chen13, Sayed14proc, Chen15transient}. {When exact gradients are employed, these strategies converge to a small area around the minimizer of the aggregate cost at a linear rate~\cite{Chen13, Yuan16}. Exact convergence requires diminishing step-sizes, resulting in sublinear rates of convergence.} A number of more recent works focusing primarily on deterministic optimization, have proposed variations yielding linear rates of convergence pursued either in the primal~\cite{Shi15, Shi15extra, DiLorenzo16, Sun16, Sun17, Xin18, Yuan19, Xin19} or dual domain~\cite{Jakovetic11, Duchi12, Tsianos12, Shi14, Jaggi14, Ling15, Jakovetic15, Seaman17, Lan18, Jakovetic19} where~\cite{Duchi12, Lan18, Xin19} allow for stochastic gradient approximations and~\cite{Jaggi14} considers empirical risk minimization problems for a linear model.

One common method for handling non-differentiable cost functions is the utilization of sub-gradient recursions, where the ordinary gradient is replaced by sub-gradients~\cite{Nedic09, Ram10, Nedic10, Duchi12, Tsianos12, Lan18}. Most often, these works assume the sub-gradients are \emph{bounded}. This condition is not satisfied in
many important cases of interest, for example, even when \(J_k(w)\) is simply quadratic in \(w\) (as happens in
mean-square-error designs) or when the \(R_k(w)\) are indicator functions used to encode constraints. Variations for
specific choices of costs functions are examined in~\cite{DiLorenzo12,DiLorenzo13,Chouvardas12,Ying12} where only the
subgradients of \(R_k(\cdot)\) are required to be bounded. The work~\cite{Ying18} generalized these conditions to allow for (sub-)gradients that are ``affine-Lipschitz'', which holds for many, but not all costs and regularizers of interest, such as indicator functions. For the case when the \(R_k(w)\) are chosen as indicator
functions in constrained problem formulations, as an alternative to projection based schemes~\cite{Nedic10, Ram10}, a distributed diffusion strategy based on the use of suitable penalty functions was proposed and studied in~\cite{Towfic15}.

Some other studies pursue distributed solutions by relying instead on the use of proximal iterations (as opposed
to sub-gradient iterations); an accessible survey on the proximal operator and its properties appears in~\cite{Boyd13}. For example, for purely deterministic costs, distributed proximal strategies are developed in~\cite{Ozdaglar12, Shi15, DiLorenzo16, Sun16, Sun17}. Stochastic variations for mean-square error costs with bounded regularizer subgradients are proposed in~\cite{Yamada13,DiLorenzo14} for single-task problems and in~\cite{Nassif15} for multi-task environments. A strategy for general stochastic costs with \emph{small}, Lipschitz continuous regularizers is studied in~\cite{Vlaski15}.

\subsection{Contributions}
\noindent The purpose of this work is to propose a general distributed strategy and a line of analysis that is applicable to a wide class of stochastic costs and non-differentiable regularizers. The first step in the solution will involve replacing each non-differentiable component, \(R_k(w)\), by a differentiable approximation \(R_k^{\delta}(w)\), parametrized by
\(\delta > 0\), such that
\begin{equation}
  \| w^o - w_{\delta}^o \|^2 \le O(\delta).
\end{equation}
The accuracy of the approximation is controlled through the smoothing parameter \( \delta \). Subsequently, we will solve for the minimizer:
\begin{equation}\label{eq:smoothproblem}
  w_{\delta}^o= \argmin_w \: \sum_{k=1}^{N} p_k \left \{ J_k(w) + R_k^{\delta}(w) \right \}
\end{equation}
Smoothing non-differentiable costs via infimal convolution~\cite{Nesterov05, Beck12, Yu13} is a popular technique in the deterministic optimization literature, and it can be used to motivate some known algorithms, such as the proximal point algorithm~\cite{Boyd13}. The technique has been mainly developed for deterministic optimization by \emph{single} stand-alone agents. In this work, we are pursue an extension in two non-trivial directions. First, we consider networked agents (rather than a single agent) working together to solve the aggregate optimization problem~\eqref{eq:wnot} (or~\eqref{eq:smoothproblem}) and, second, the risk functions involved are a combination of stochastic costs defined as the expectations of certain loss functions and deterministic regularizers. Moreover, the probability distribution of the data is assumed unknown and, therefore, the aggregate risks themselves are not known but can only be approximated. The challenge is to devise a distributed strategy that is able to converge to the desired Pareto solution despite these difficulties.

We note that an alternative smoothing procedure by means of adding small stochastic perturbations is considered in~\cite{Duchi12smoothing} and extended to decentralized stochastic optimization in~\cite{Scaman18}, requiring bounded subgradients. In contrast, our focus is on smooth stochastic risks regularized by non-smooth, deterministic risks. Splitting the smooth stochastic part from the non-differentiable deterministic risk, and smoothing only the deterministic risk via a deterministic procedure will allow us to only require looser bounds on both components.

In the next sections we will explain how to construct the smooth approximation, \(R_{k}^{\delta}(w)\), by
appealing to conjugate functions and will show that the distance \( \|w^o - w_{\delta}^o\| \) can be made arbitrarily small for \(\delta \to 0\). We then present an algorithm to solve for the minimizer of~\eqref{eq:smoothproblem} in a distributed manner and derive bounds on its performance. The analysis in future sections will rely on the following common assumptions~\cite{Sayed14proc,Sayed14, Chen15transient}:
\begin{assumption}[Lipschitz gradients]\label{as:lipschitz_gradient_general_regularizer}
  For each \(k\), the gradient \(\nabla J_k(\cdot)\) is Lipschitz, namely, there exists \( \lambda_U \geq 0 \) such that for any \(x,y \in \mathds{R}^{M}\):
  \begin{equation}\label{eq:gradbound}
    \|\nabla J_k(x) - \nabla J_k(y)\| \le \lambda_{U} \|x-y\|
  \end{equation}\hfill\IEEEQED
\end{assumption}
\begin{assumption}[Strong Convexity]\label{as:strongcon}
  The weighted aggregate of the differentiable risks is strongly convex, namely, there exists \( \lambda_L \geq 0 \) such that for any \(x,y \in \mathds{R}^{M}\):
  \begin{equation}\label{eq:jkstrong}
    {\left(x-y\right)}^{\T} \cdot \sum_{k=1}^{N} p_k\left(\nabla_w J_k(x) - \nabla_w J_k(y)\right) \ge \lambda_L \|x-y\|^2
  \end{equation}\hfill\IEEEQED
\end{assumption}
\begin{assumption}[Regularizers]\label{as:reg}
  For each \(k\), \(R_k(\cdot)\) is closed convex. In other words, \( R_k(\cdot) \) is convex and\\ \( \left \{ w \in \mathrm{dom}\,R_k(\cdot)\,|\,R_k(w) \le x \right \} \) is a closed set for all \( x \). \hfill\IEEEQED
\end{assumption}

\section{Algorithm Formulation}
\subsection{Construction of Smooth Approximation}
\noindent To begin with, following the works~\cite{Nesterov05, Beck12}, we explain how smoothing of the regularizers is performed. Thus, recall  that the conjugate function,
denoted by \( R_k^{\star}(w) \), of a regularizer \( R_k(w) \) is defined as
\begin{equation}\label{eq:dualrep}
  R_k^{\star}(w) \triangleq \sup_{u \in \dom R_k} \left \{ w^{\T} u - R_k(u) \right \}.
\end{equation}
A useful property of conjugate functions is that \( R_{k}^{\star}(w) \) is always closed convex regardless of whether
\( R_k(w) \) is convex or not.
\begin{definition}[Proximity function~\cite{Nesterov05}]\label{def:distance}
  A proximity function \( d(\cdot) \) for a closed convex set \( C \) is a continuous, strongly-convex function with
  \( C \subseteq \dom d(\cdot) \). We center and normalize
  the function so that
  \begin{equation}\label{eq:distinf}
    \min_{w \in C} \: d(w) = 0
  \end{equation}
  and
  \begin{equation}
    \argmin_{w \in C} \: d(w) = 0
  \end{equation}
  which exists and is unique, since \( d(w) \) is strongly-convex. Furthermore, the proximity function is scaled
  to satisfy the following normalization (which means that its strong-convexity constant is set to one):
  \begin{equation}
     d(w) \ge \frac{1}{2} \|w\|^2.
  \end{equation}\hfill\IEEEQED
\end{definition}
\begin{definition}[Smooth approximation]
  We choose a proximity function over \( C = \dom R_k^{\star}(w) \) and define the smooth approximation of
  \( R_k(\cdot) \) as:
  \begin{align}
    R_k^{\delta}(w) \triangleq& \max_{u \in \dom R_k^{\star}} \: \left \{ w^{\T} u - R_k^{\star}(u) - \delta\cdot d(u) \right \} \notag \\
    =& {\left( R_k^{\star} + \delta\cdot d\right)}^{\star}\left(w\right) \label{eq:smoothshort}
  \end{align}\hfill\IEEEQED
\end{definition}
\noindent
The maximum in~\eqref{eq:smoothshort} is attained for all \( w \) since \( R_k^{\star}(u) + \delta\cdot d(u) \) is
strongly convex. Thus, observe that the smooth approximation for \( R_k(w) \), which we are denoting by \(
R_{k}^{\delta}(w) \), is obtained by first perturbing the conjugate function \( R_k^{\star}(u) \) by \( \delta\cdot d(u) \)
and then conjugating the result again. The perturbation makes the sum \( R_{k}^{\star}(u)+\delta\cdot d(u) \) a
strongly-convex function. The motivation behind this construction is the fact that the conjugate of a strongly-convex
function is \emph{differentiable} everywhere and, therefore, \( R_{k}^{\delta}(w) \) is differentiable everywhere. This
intuition is formalized in the following known theorem~\cite{Nesterov05}, preceded by an elementary lemma~{\cite{Rockafellar70}}.
\begin{lemma}[Conjugate subgradients~{\cite{Rockafellar70}}]\label{lem:subgradients}
  If \( G(\cdot) \) is some closed and convex function, the subgradients of \( G(\cdot) \) and its conjugate \( G^{\star}(\cdot) \)
  are related as:
  \begin{equation}
    v \in \partial G(w) \longleftrightarrow w \in \partial G^{\star}(v)
  \end{equation}
\end{lemma}
\begin{IEEEproof}
  The lemma is from~{\cite{Rockafellar70}}.\ifarx For reference, the proof is repeated in Appendix~\ref{ap:subgradients}.\fi
\end{IEEEproof}
\begin{theorem}[Gradient of smooth approximation~\cite{Nesterov05}]\label{th:smoothgradient}
  Any \( R_k^{\delta}(w) \) constructed according to~\eqref{eq:smoothshort} is differentiable with gradient vector
  \begin{equation}\label{eq:gradient}
    \nabla R_k^{\delta}(w) = \underset{u \in \dom R_k^{\star}}{\arg\max} \: \left \{
      w^{\T} u - R_k^{\star}(u) - \delta\cdot d(u) \right \}.
  \end{equation}
  Furthermore, the gradient is co-coercive, i.e., it satisfies:
  \begin{equation}\label{eq:rkcoco}
    {\left(x-y\right)}^{\T}\left(\nabla R_k^{\delta}(x) - \nabla R_k^{\delta}(y)\right) \ge \delta \|\nabla R_k^{\delta}(x) - \nabla R_k^{\delta}(y)\|^2
  \end{equation}
  By Cauchy-Schwarz, this implies Lipschitz continuity, i.e.,
  \begin{equation}\label{eq:rklip}
    \|\nabla R_k^{\delta}(x) - \nabla R_k^{\delta}(y)\| \le \frac{1}{\delta} \|x-y\|.
  \end{equation}
\end{theorem}
\begin{IEEEproof}
  The theorem is from~\cite{Nesterov05}. For reference, the proof is repeated in Appendix~\ref{ap:smoothgradient}.
\end{IEEEproof}
The feasibility of stochastic-gradient algorithms for the minimization of~\eqref{eq:smoothproblem}
hinges on the assumption that~\eqref{eq:gradient} can be evaluated in closed form or at least easily. Fortunately, this
is the case for a large class of regularizers of interest --- see~\cite{Combettes11} for an overview of closed form
solutions in the special case \( d(\cdot)=\frac{1}{2}{\| \cdot \|}^2 \) and~\cite{Nesterov05, Beck12} for other distance choices. {For example, for every function where the proximal operator~\cite{Boyd13}:
\begin{equation}
  R_k^{\delta}(w) = \underset{u}{\min} \: \left( R_k(w) + \frac{1}{2\delta} \|w-u\|^2 \right)
\end{equation}
can be evaluated in closed form, we can let \( d(\cdot) \triangleq \frac{1}{2} \| \cdot\|^2 \) and obtain~\cite{Boyd13}:
\begin{equation}
  \nabla R_k^{\delta}(w) = \frac{1}{\delta} \left( w - \mathrm{prox}_{\delta R_k}(w) \right)
\end{equation}
Depending on the regularizers \( R_k(\cdot) \), other proximity functions may be more appropriate~\cite{Nesterov05}}. We point out that the smooth approximation~\eqref{eq:smoothshort} can equivalently be written as~\cite{Beck12}:
\begin{align}\label{eq:inf_convolution}
  R_k^{\delta}(w) &= \min_{u \in \dom R_k} \: \left \{ R_k(u) + \delta\cdot d^{\star}\left(\frac{w - u}{\delta}\right) \right \}
\end{align}\ifarx
To verify this, observe that
\begin{align}
  R_k^{\delta}(w) &= \min_{u \in \dom R_k} \: \left \{ R_k(u) + \delta\cdot \sup_z \left \{ z^{\T} \left( \frac{w-u}{\delta} \right) - d\left(z\right) \right \} \right \} \notag \\
  &= \min_{u \in \dom R_k} \: \left \{ R_k(u) + \sup_z \left \{ z^{\T} \left( w-u \right) - \delta \cdot d\left(z\right) \right \} \right \} \notag \\
  &= \sup_z \: \left \{ \inf_u \left \{- z^{\T} u + R_k(u) \right \} + z^{\T} w -  \delta \cdot d\left(z\right) \right \} \notag \\
  &= \sup_z \: \left \{ - \sup_u \left \{R_k(u) - z^{\T} u \right \} + z^{\T} w -  \delta \cdot d\left(z\right) \right \} \notag \\
  &= \max_z \: \left \{ z^{\T} w - R_k^{\star}(z) -  \delta \cdot d\left(z\right) \right \}
\end{align}\fi
Expression~\eqref{eq:inf_convolution} is known as the infimal convolution.

\subsection{Accuracy of the Smooth Approximation}
\noindent Replacing the original optimization problem~\eqref{eq:wnot} by the smoothed cost~\eqref{eq:smoothproblem} naturally results in a bias, since the new minimizer \( w_{\delta}^o \) will generally be different from the original minimizer \( w^o \). This bias, when not properly controlled, can degrade the performance of the algorithm. For this reason, a number of works have examined the smoothing bias introduced through conjugate smoothing under various conditions on the cost functions. In the centralized setting, when \( N = 1 \), it has been established that \( R_k^{\delta}(w) \to R_k(w) \) both pointwise and epigraphically, which implies \( w_{\delta}^o \to w^o \) as \( \delta \to 0 \)~\cite{Planiden16}, while~\cite{Bauschke08} showed a sum of costs \( \sum_{k=1}^N p_k R_k(w) \), when smoothed individually, will continue to converge epigraphically. While encouraging, these results do not guarantee a rate at which \( w_{\delta}^o \to w^o \), complicating the choice of the smoothing parameter \( \delta \). Pointwise convergence has been strengthened to uniform convergence, i.e., \( \left|R_k(w) - R_k^{\delta}(w)\right| \le O(\delta) \) for costs with bounded subgradients for \( N =1 \)~\cite{Nesterov05, Beck12} and for a collection of costs, each with bounded subgradients in~\cite{Yu13}.

We present here a variation of these results by restricting ourselves to \emph{strongly}-convex costs, but allowing for regularizers with \emph{unbounded} sub-gradients and establishing \( \|w^o - w_{\delta}^o\|^2 \le O(\delta) \) rather than simply \( w_{\delta}^o \to w^o \).
\begin{theorem}[Accuracy of smooth approximation]\label{th:smoothacc}
  The bias introduced by smoothing the original problem diminishes linearly with \( \delta \), i.e.,
  \begin{equation}\label{eq:minimizer_bias}
    \|w^o - w_{\delta}^o\|^2 \le \delta \cdot \frac{2}{\lambda_L} \sum_{k=1}^{N} p_k d\left( r_k^o \right) = O(\delta)
  \end{equation}
  where \( r_k^o \in \partial R_k(w^o) \) such that
  \begin{equation}
    \sum_{k=1}^N p_k \left \{ \nabla J_k(w^o) + r_k^o \right \} = 0.
  \end{equation}
  This collection of \( \left \{ r_k^o \right \}\) is guaranteed to exist, since \( w^o \triangleq \argmin \sum_{k=1}^N p_k \left \{ J_k(w) + R_k(w) \right \} \).
\end{theorem}
\begin{IEEEproof} Appendix~\ref{ap:smoothacc}. \end{IEEEproof}

\subsection{Regularized Diffusion Strategy}
\noindent Now that we have established a method for constructing a differentiable approximation for each regularizer, we can solve
for the minimizer of~\eqref{eq:smoothproblem} by resorting to the following (adapt-then-combine form of the) diffusion
strategy~\cite{Sayed14proc,Sayed14,Chen15transient}:
\begin{align}
  \boldsymbol{\phi}_{k,i} &= \w_{k,i-1} - \mu \widehat{\nabla J}_{k}(\w_{k,i-1}) - \mu \nabla R_{k}^{\delta}(\w_{k,i-1})\label{eq:org1}\\
  \w_{k,i} &= \sum_{\ell=1}^{N} a_{\ell k}\boldsymbol{\phi}_{\ell,i} \label{eq:org2}
\end{align}
where \( \mu>0 \) is a small step-size parameter {and \( a_{\ell k} \) are the entries of a combination matrix \( A \) with Perron eigenvector \( p \), i.e. \( A p = p \)}. In this implementation, each agent \( k \) first performs the
stochastic-gradient update~\eqref{eq:org1}, starting from its existing iterate value \( \w_{k,i-1} \), and obtains an
intermediate iterate \( \boldsymbol{\phi}_{k,i} \). Subsequently, agent \( k \) consults with its neighbors and combines their
intermediate iterates into \( \w_{k,i} \) according to~\eqref{eq:org2}. Motivated by the construction in~\cite{Towfic15}, we
can refine~\eqref{eq:org1}--\eqref{eq:org2} further as follows. We first introduce an auxiliary variable
\( \boldsymbol{\psi}_{k,i} \) and rewrite~\eqref{eq:org1} in the equivalent form:
\begin{align}
  \boldsymbol{\phi}_{k,i} &= \w_{k,i-1} - \mu \widehat{\nabla J}_{k}(\w_{k,i-1})\\
  \boldsymbol{\psi}_{k,i} &= \boldsymbol{\phi}_{k,i} - \mu \nabla R_{k}^{\delta}(\w_{k,i-1}) \label{eq:intermediateStep}\\
  \w_{k,i} &= \sum_{\ell=1}^{N} a_{\ell k}\boldsymbol{\psi}_{\ell,i}
\end{align}
We can now appeal to an incremental-type argument~\cite{Bertsekas97incremental,Bertsekas97parallel} by noting that it is
reasonable to expect \( \boldsymbol{\phi}_{k,i} \) to be an improved estimate for \( w_{\delta}^o \) compared to \(
\w_{k,i-1} \). Therefore, we replace \( \w_{k,i-1} \) in~\eqref{eq:intermediateStep} by \( \boldsymbol{\phi}_{k,i} \)
and arrive at the following regularized diffusion implementation.
\begin{algorithm}
  \SetAlgoLined%
  \begin{align}
    \boldsymbol{\phi}_{k,i} &= \w_{k,i-1} - \mu \widehat{\nabla J}_{k}(\w_{k,i-1})\label{eq:adapt1}\\
    \boldsymbol{\psi}_{k,i} &= \boldsymbol{\phi}_{k,i} - \mu \nabla R_{k}^{\delta}(\boldsymbol{\phi}_{k,i})\label{eq:adapt2}\\
    \w_{k,i} &= \sum_{\ell=1}^{N} a_{\ell k} \boldsymbol{\psi}_{\ell,i}\label{eq:combine}
  \end{align}
  \caption{Regularized Diffusion Strategy}
\end{algorithm}
\begin{example}[Proximal Diffusion Learning]
  Choosing \( d(w)=\frac{1}{2}\|w\|^2 \) turns the smooth approximation~\eqref{eq:smoothshort} into
  \begin{equation}
    R_k^{\delta}(w) = {\left( R_k^{\star}(w) + \frac{\delta}{2} \|w\|^2\right)}^{\star}
  \end{equation}
  which is the well-known Moreau envelope~\cite{Boyd13}. It can be rewritten equivalently as
  \begin{equation}
    R_k^{\delta}(w) = \underset{u}{\min} \: \left( R_k(w) + \frac{1}{2\delta} \|w-u\|^2 \right)
  \end{equation}
  where the minimizing argument is identified as the proximal operator:
  \begin{equation}
    \mathrm{prox}_{\delta R_k}(w) = \underset{u}{\arg\min} \: \left( R_k(w) + \frac{1}{2\delta} \|w-u\|^2 \right)
  \end{equation}
  For many costs \( R_k(w) \), the proximal operator can be evaluated in closed form. The gradient of the Moreau envelope
  can also be written as
  \begin{equation}
    \nabla R_k^{\delta}(w) = \frac{1}{\delta} \left( w - \mathrm{prox}_{\delta R_k}(w) \right).
  \end{equation}
  This allows us to rewrite iterations~\eqref{eq:adapt1}--\eqref{eq:combine} as
  \begin{align}
    \boldsymbol{\phi}_{k,i} &= \w_{k,i-1} - \mu \widehat{\nabla J}_{k}(\w_{k,i-1})\\
    \boldsymbol{\psi}_{k,i} &= \left(1-\frac{\mu}{\delta}\right) \boldsymbol{\phi}_{k,i} + \frac{\mu}{\delta}\mathrm{prox}_{\delta R_k} (\boldsymbol{\phi}_{k,i})\\
    \w_{k,i} &= \sum_{\ell=1}^{N} a_{\ell k}\boldsymbol{\psi}_{\ell,i}
  \end{align}
  which is a damped variation of the proximal diffusion algorithm studied in~\cite{Vlaski15} under
  the stronger assumption of small Lipschitz continuous regularizers.\hfill\IEEEQED%
\end{example}

\section{Convergence Analysis}

\subsection{Centralized Recursion}
\noindent We now examine the convergence properties of the diffusion strategy~\eqref{eq:adapt1}--\eqref{eq:combine}. To do so, and motivated by the approach introduced in~\cite{Chen15transient}, it
is useful to introduce the following centralized recursion to serve as a frame of reference:
\begin{equation}\label{eq:centralized}
  w_i = w_{i-1} - \mu \sum_{k=1}^{N} p_k \nabla J_k(w_{i-1}) - \mu \sum_{k=1}^{N} p_k \nabla R_k^{\delta}(w_{i-1})
\end{equation}
This recursion amounts to a gradient-descent iteration applied to the smoothed aggregate cost
in~\eqref{eq:smoothproblem} under the assumption that the risk functions (and therefore their gradients) are known. For convenience of presentation, we
introduce the central operator \( T_c(x):\mathds{R}^{M}\rightarrow\mathds{R}^{M} \) defined as follows:
\begin{align}\label{eq:Tc}
  T_c(x)&\triangleq x - \mu \sum_{k=1}^{N} p_k \nabla J_k(x) - \mu \sum_{k=1}^{N} p_k \nabla R_k^{\delta}(x)
\end{align}
so that the reference recursion~\eqref{eq:centralized} becomes \( w_i=T_c(w_{i-1}) \).
\begin{lemma}[Contraction mapping]\label{lem:centcontr}
  Assume \( \mu \le 2\delta \). Then, the centralized recursion operator~\eqref{eq:Tc} satisfies
  \begin{equation}
    \|T_c(x)-T_c(y)\| \le \gamma_c \|x-y\|
  \end{equation}
  where \( \gamma_c>0 \) can be made strictly less than one by selecting sufficiently small \( \mu \) and is given by:
  \begin{equation}
    \gamma_c=1- \mu \lambda_L + \mu^2
    \left(\frac{\lambda_U^2}{2-\frac{\mu}{\delta}}\right).
  \end{equation}
  From Banach's fixed point theorem~\cite{Kreyszig89} and~\eqref{eq:Tc}, we conclude that for sufficiently small \( \mu \),
  \( w_i=T_c(w_{i-1}) \) converges exponentially to the unique fixed-point \( w^o_{\delta} \), the minimizer of~\eqref{eq:smoothproblem}.
\end{lemma}
\begin{IEEEproof} Appendix~\ref{ap:centcontr}.\end{IEEEproof}

\subsection{Network Basis Transformation}
\noindent We are now ready to examine the behavior of the diffusion strategy~\eqref{eq:adapt1}--\eqref{eq:combine}, which employs stochastic gradients. We begin by
introducing the following extended vectors and matrices, which collect quantities of interest from across all agents in
the network:
\begin{align}
  \bcw_{i} &\triangleq \mathrm{col} \left \{ \w_{1,i}, \ldots, \w_{N,i} \right \} \\
  \mathcal{A} &\triangleq A \otimes I_M\\
  {g}(\bcw_{i}) &\triangleq \mathrm{col} \left \{ {\nabla_w J_1}(\w_{1,i}), \ldots ,{\nabla_w J}_N( \w_{N,i}) \right \} \\
  \widehat{g}(\bcw_{i}) &\triangleq \mathrm{col} \left \{ \widehat{\nabla_w J}_1(\w_{1,i}), \ldots , \widehat{\nabla_w J}_N(\w_{N,i}) \right \} \\
  r(\bcw_{i}) &\triangleq \mathrm{col} \left \{ \nabla_w R_1^{\delta} (\w_{1,i}), \ldots ,\nabla_w R_N^{\delta}(\w_{N,i}) \right \} \\
  q(\bcw_{i}) &\triangleq r(\bcw_{i}-\mu {g} (\bcw_{i}))\\
  \widehat{q}(\bcw_{i}) &\triangleq r(\bcw_{i}- \mu \widehat{g}(\bcw_{i}))
\end{align}
Using these definitions, iterations~\eqref{eq:adapt1}--\eqref{eq:combine} show that the network vector
\( \bcw_i \) evolves according to the following dynamics:
\begin{equation}\label{eq:fullrecursion}
  \bcw_i = \mathcal{A}^{\T} \bcw_{i-1} - \mu \mathcal{A}^{\T}
  \left(\widehat{g}(\bcw_{i-1})+ \widehat{q}(\bcw_{i-1})\right)
\end{equation}
By construction, the combination matrix \( A \) is left-stochastic and primitive and hence admits a Jordan decomposition of
the form \( A = V_{\epsilon} J V_{\epsilon}^{-1} \) with~\cite{Sayed14,Chen15transient}:
\begin{equation}\label{eq:jordan_general_regularizer}
  V_{\epsilon} = \left[ \begin{array}{cc} p & V_R \end{array} \right],
  \ \ J = \left[ \begin{array}{cc} 1 & 0\\0 & J_{\epsilon} \end{array}\right],
  \ \ V_{\epsilon}^{-1} = \left[ \begin{array}{c} \mathds{1}^{\T} \\ \vphantom{O^{O^{O^O}}} V_L^{\T} \end{array}\right]
\end{equation}
where \( J_{\epsilon} \) is a block Jordan matrix with the eigenvalues \( \lambda_2(A) \) through \( \lambda_N(A) \) on
the diagonal and \( \epsilon \) on the first lower sub-diagonal. The extended matrix \( \mathcal{A} \) then satisfies \(
\mathcal{A} = \mathcal{V}_{\epsilon} \mathcal{J} \mathcal{V}_{\epsilon}^{-1} \) with \( \mathcal{V}_{\epsilon} =
{V}_{\epsilon} \otimes I_N \), \( \mathcal{J} = J \otimes I_N \), \( \mathcal{V}_{\epsilon}^{-1} = {V}_{\epsilon}^{-1}
\otimes I_N \). Multiplying both sides of~\eqref{eq:fullrecursion} by \( \mathcal{V}_{\epsilon}^{\T} \) and introducing
the transformed iterate vector \( \bcw'_i \triangleq \mathcal{V}_{\epsilon}^{\T} \bcw_i \), we obtain
\begin{equation}\label{eq:transformedrecursion}
  \bcw'_i = \mathcal{J}^{\T} \bcw'_{i-1} - \mu \mathcal{J}^{\T}
  \mathcal{V}_{\epsilon}^{\T} \left(\widehat{g}(\bcw_{i-1})+\widehat{q}(\bcw_{i-1}) \right)
\end{equation}
Following~\cite{Chen15transient, Sayed14}, we can exploit the structure of the decomposition~\eqref{eq:jordan_general_regularizer} to provide further insight into this transformed recursion. Let \( {\bcw'_i} = \mathrm{col} \left \{ \w_{c,i}, \bcw_{e,i}\right \} \), where \( \w_{c,i} \in
\mathds{R}^{N \times 1} \) and \( \bcw_{e,i} \in \mathds{R}^{(N-1)M \times 1} \). Then,
recursion~\eqref{eq:transformedrecursion} can be decomposed as
\begin{align}
  \w_{c,i} =&\: \w_{c,i-1} - \mu \left( p^{\T} \otimes I_N\right) \left(\widehat{g}(\bcw_{i-1})+\widehat{q} (\bcw_{i-1})\right) \label{eq:centrecursion}\\
  \bcw_{e,i} =&\: \mathcal{J}_{\epsilon}^{\T} \bcw_{e,i-1} - \mu \mathcal{J}_{\epsilon}^{\T} \mathcal{V}_{R}^{\T} \left(\widehat{g}(\bcw_{i-1}) + \widehat{q}(\bcw_{i-1})\right) \label{eq:errorrecursion}
\end{align}
Note from \( \bcw'_i = \mathcal{V}_{\epsilon}^{\T} \bcw_i \), that~\cite{Chen15transient}:
\begin{equation}
  \w_{c,i} = \left( p^{\T} \otimes I_M \right) \bcw_i = \sum_{k=1}^{N} p_k
  \w_{k,i}
\end{equation}
Hence, \( \w_{c,i} \) is the weighted centroid vector of all iterates \( \w_{k,i} \) across the network. From \( \bcw_i =
{\left(\mathcal{V}_{\epsilon}^{-1}\right)}^{\T} \bcw'_i \) on the other hand, one obtains~\cite{Chen15transient}:
\begin{equation}
  \bcw_i = \mathds{1} \otimes \w_{c,i} + \mathcal{V}_L \bcw_{e,i}
\end{equation}
so that \( \bcw_{e,i} \) can be interpreted as the deviation of individual estimates from the weighted centroid vector
\( \w_{c,i} \) across the network.

We examine the centroid recursion~\eqref{eq:centrecursion} in greater detail. Thus, note that
\begin{align}
  &\: \w_{c,i} \notag \\
  =&\: \w_{c,i-1} - \mu \left( p^{\T} \otimes I_M\right)
  \left(\widehat{g}(\bcw_{i-1})+\widehat{q}(\bcw_{i-1}) \right)\notag \\
  =&\: \w_{c,i-1} - \mu \left( p^{\T} \otimes I_M\right) \left({g}(\mathds{1} \otimes \w_{c,i-1}) +{r}(\mathds{1} \otimes \w_{c,i-1})\right) \notag \\
  &- \mu \left( p^{\T} \otimes I_M\right) \big({g}(\bcw_{i-1}) + {q}(\bcw_{i-1}) \notag \\
  &\: \ \ \ \ \ \ \ \ \ \ \ \ \ \ \ \ \ \ \ \  - {g}(\mathds{1} \otimes \w_{c,i-1}) - {q}(\mathds{1} \otimes \w_{c,i-1})\big) \notag \\
  &- \mu \left( p^{\T} \otimes I_M\right)\big(\widehat{g} (\bcw_{i-1}) + \widehat{q} (\bcw_{i-1}) \notag \\
  &\: \ \ \ \ \ \ \ \ \ \ \ \ \ \ \ \ \ \ \ \ -{g}( \bcw_{i-1}) - {q}( \bcw_{i-1})\big) \notag \\
  &- \mu \left( p^{\T} \otimes I_M\right) \left({q}(\bcw_{i-1})-{r}( \bcw_{i-1})\right)\notag \\
  =&\: T_c(\w_{c,i-1}) - \mu \left( p^{\T} \otimes I_M\right)\big( \boldsymbol{t}_{i-1} + \s_i + \boldsymbol{u}_{i-1} \big) \label{eq:pertcent}
\end{align}
where we replaced
\begin{align}
  \w&_{c,i-1} - \mu \left( p^{\T} \otimes I_M\right) \left({g}(\mathds{1} \otimes \w_{c,i-1})+{r}(\mathds{1} \otimes \w_{c,i-1})\right)\notag \\
  \stackrel{\phantom{(36)}}{=}& \w_{c,i-1} - \mu \sum_{k=1}^{N} p_k \nabla J_k(\w_{c,i-1}) - \mu \sum_{k=1}^{N} p_k \nabla R_k^{\delta}(\w_{c,i-1})\notag \\
  \stackrel{(\ref{eq:Tc})}{=}& T_c(\w_{c,i-1})
\end{align}
and introduced the perturbation terms:
\begin{align}
  \boldsymbol{t}_{i-1} =&\: {g}(\bcw_{i-1}) + {q}(\bcw_{i-1}) - {g}(\mathds{1} \otimes \w_{c,i-1}) - {q}(\mathds{1} \otimes \w_{c,i-1})\\
  \s_{i} =&\: \widehat{g} (\bcw_{i-1}) + \widehat{q} (\bcw_{i-1})-{g}( \bcw_{i-1}) - {q}( \bcw_{i-1})\\
  \boldsymbol{u}_{i-1} =&\: {q}(\bcw_{i-1})-{r}( \bcw_{i-1})
\end{align}
It follows from~\eqref{eq:pertcent} that the centroid recursion is a perturbed version of the central recursion
introduced earlier in~\eqref{eq:Tc}. The perturbation arising from disagreement across agents in the network is captured
in \( \boldsymbol{t}_{i-1} \), while stochastic perturbations due to instantaneous gradient approximations is captured
in \( \s_i \). The incremental implementation causes \( \boldsymbol{u}_{i-1} \). It is therefore reasonable to expect
that \( \w_{c,i} \) will evolve close to the central variable \( w_i \) from~\eqref{eq:centralized}, which was already
shown to converge to \( w_{\delta}^o \) in Lemma~\ref{lem:centcontr}. To formalize this intuition, we define \(
\widetilde{\w}_{c,i-1} = w_{\delta}^o - \w_{c,i-1} \). Since \( w_{\delta}^o \) is a fixed-point of \( T_c(\cdot) \),
i.e., \( w_{\delta}^o = T_c(w_{\delta}^o) \), the error \( \widetilde{\w}_{c,i-1} \) satisfies the recursion
\begin{align}
  \widetilde{\w}_{c,i-1}=&\: T_c(w_{\delta}^o) - T_c(\w_{c,i-1}) \notag \\
  &\: + \mu \left( p^{\T} \otimes I_M\right)\big( \boldsymbol{t}_{i-1} + \s_i + \boldsymbol{u}_{i-1} \big) \label{eq:perterror}
\end{align}
With the same perturbation terms, expression~\eqref{eq:errorrecursion} turns
into
\begin{align}
  \bcw_{e,i} =&\: \mathcal{J}_{\epsilon}^{\T} \bcw_{e,i-1} - \mu \mathcal{J}_{\epsilon}^{\T} \mathcal{V}_{R}^{\T}\big( \boldsymbol{t}_{i-1} + \s_i + \boldsymbol{u}_{i-1} \notag \\
  &\: -{g}(\mathds{1} \otimes \w_{c,i-1})-{r}(\mathds{1} \otimes \w_{c,i-1})\big)
\end{align}
We employ the following common assumption on the perturbations caused by the gradient noise~\cite{Sayed14proc, Sayed14,
Chen15transient}.
\begin{assumption}[Gradient noise process]\label{eq:gradientnoise_general_regularizer}
  For each \( k \), the gradient noise process is defined as
  \begin{equation}
    \s_{k,i}(\w_{k,i-1}) = \widehat{\nabla J}_k(\w_{k,i-1}) - \nabla J_k(\w_{k,i-1})
  \end{equation}
  and satisfies
  \begin{subequations}
    \begin{align}
      \E \left[ \s_{k,i}(\w_{k,i-1}) | \boldsymbol{\mathcal{F}}_{i-1} \right] &= 0\\
      \E \left[ \|\s_{k,i}(\w_{k,i-1})\|^2 | \boldsymbol{\mathcal{F}}_{i-1} \right] &\le \beta^2 \|\w_{k,i-1}\|^2 + \sigma^2 \label{eq:noisebound}
    \end{align}
  \end{subequations}
  for some non-negative constants \( \{\beta^2,\sigma^2\} \), and where \( \boldsymbol{\mathcal{F}}_{i-1} \) denotes the filtration
  generated by the random processes \( \{\w_{\ell,j}\} \)  for all \( \ell=1,2,\ldots,N \) and \( j \le i-1 \), i.e.,
  \( \boldsymbol{\mathcal{F}}_{i-1} \) represents the information that is available about the random processes
  \( \{\w_{\ell,j}\} \) up to time \( i-1 \).\hfill\IEEEQED%
\end{assumption}
\noindent For a block-vector \( \boldsymbol{x} \in \mathds{R}^{MN\times1} \) consisting of \( N \) blocks of size \(
M\times 1 \), let \( P[\boldsymbol{x}]=\mathrm{col}\left \{ \E \|\boldsymbol{x}_1\|^2, \dots, \E \|\boldsymbol{x}_N\|^2
\right \} \in \mathds{R}^{N\times 1} \)~\cite{Chen15transient}. Note that \( \mathds{1}^{\T} P[\boldsymbol{x}] = \E \|\boldsymbol{x}\|^2 \).
Furthermore, let \( v_{L,k} \) denote the \( k \)-th row of \( V_L \) and let \( \nu = \max_k \: \|v_{L,k} \otimes I_M\| \),
which is independent of \( \mu \) and \( \delta \).
\begin{lemma}[Bounds on perturbation terms]\label{lem:pertbounds}
  The perturbation terms in~\eqref{eq:perterror} satisfy the following bounds:
  \begin{align}
    &P[\boldsymbol{t}_{i-1}] \preceq \left(2 \lambda_U^2 + 4 \frac{1+\mu^2}{\delta^2} \right) \nu^2 \mathds{1}\mathds{1}^{\T} P[\bcw_{e,i-1}]\label{eq:boundt}\\
    &P[\boldsymbol{u}_{i-1}] \preceq \frac{\mu^2}{\delta^2} \big(3\lambda_U^2 \nu^2 \mathds{1}\mathds{1}^{\T}P[\bcw_{e,i-1}] + 3\lambda_U^2 P[\mathds{1}\otimes \widetilde{\w}_{c,i-1}] \notag \\
    &\ \ \ \ \ \ \ \  \ \ \ \ \   + 3 P[g(\mathds{1} \otimes w_{\delta}^o)]\big)\label{eq:boundu}\\
    &P[\s_{i} -  \E \s_i] \preceq 3\beta^2 P[\mathds{1}\otimes \widetilde{\w}_{c,i-1}] + 3 \beta^2 \nu^2 \mathds{1}\mathds{1}^{\T}P[\bcw_{e,i-1}] \notag \\
    & \ \ \ \ \ \ \ \ \ \ \ \ \ \  \ \ \ \   + 3\beta^2 P[\mathds{1}\otimes w_{\delta}^o] + \sigma^2 \mathds{1}\label{eq:bounds}\\
    &P[\E \s_i] \preceq 3\beta^2 \frac{\mu^2}{\delta^2} P[\mathds{1}\otimes \widetilde{\w}_{c,i-1}]  + 3 \beta^2 \frac{\mu^2}{\delta^2} \nu^2 \mathds{1}\mathds{1}^{\T}P[\bcw_{e,i-1}]\notag \\
    &  \ \ \ \ \ \ \ \ \ \ \ \ + 3\beta^2 \frac{\mu^2}{\delta^2}P[\mathds{1}\otimes w_{\delta}^o] + \frac{\mu^2}{\delta^2}\sigma^2 \mathds{1}\label{eq:bounde}\\
    &P[g(\mathds{1}\otimes \w_{c,i-1})] \preceq 2 \lambda_U^2 P[\mathds{1}\otimes \widetilde{\w}_{c,i-1}] + 2 P[g(\mathds{1}\otimes w_{\delta}^o)]\label{eq:boundgc}\\
    &P[r(\mathds{1}\otimes \w_{c,i-1})] \preceq \frac{2}{\delta^2} P[\mathds{1}\otimes \widetilde{\w}_{c,i-1}] + 2 P[r(\mathds{1}\otimes w_{\delta}^o)]\label{eq:boundrc}
  \end{align}
\end{lemma}
\begin{IEEEproof} Appendix~\ref{ap:pertbounds}.\end{IEEEproof}

\subsection{Mean-Square-Error Bounds}
\noindent Using the bounds on the perturbation terms obtained in Lemma~\ref{lem:pertbounds}, we can formulate a recursive bound on
the mean-square error.
\begin{lemma}[Mean-Square-Error Recursion]\label{lem:errorrec}
  \allowdisplaybreaks%
  The variances of \( \widetilde{\w}_{c,i} \) and \( \bcw_{e,i} \) are coupled and recursively bounded as
  \begin{equation}
    \begin{bmatrix}
      \E {\|\widetilde{\w}_{c,i}\|}^2\\
      \E {\|\bcw_{e,i}\|}^2
    \end{bmatrix}
    \preceq \Gamma \begin{bmatrix}
      \E \|\widetilde{\w}_{c,i-1}\|^2\\
      \E \|\bcw_{e,i-1}\|^2
    \end{bmatrix} + \begin{bmatrix}
      \frac{\mu^3}{\delta^2}b_1 + \frac{\mu^3}{\delta^2}b_2 + \mu^2 b_3\\
      \frac{\mu^2}{\delta^2}b_4 + \mu^2 b_5 + \frac{\mu^4}{\delta^2} b_6
    \end{bmatrix}
  \end{equation}
  where
  \begin{align}
    \Gamma=&\:
    \begin{bmatrix} \gamma_c + \frac{\mu^3}{\delta^2} h_1 + \mu^2 h_2 & \frac{\mu}{\delta^2} h_3 + \mu h_4 + \frac{\mu^3}{\delta^2} h_5 + \mu^2 h_6\\
      \frac{\mu^2}{\delta^2} h_7 + \mu^2 h_8 + \frac{\mu^4}{\delta^2} h_9 & \|\mathcal{J}_{\epsilon}\|+\frac{\mu^2}{\delta^2} h_{10} + \mu^2 h_{11} + \frac{\mu^4}{\delta^2} h_{12}
    \end{bmatrix}\label{eq:driving_matrix}\\
    \gamma_c\triangleq&\:1- \mu \lambda_L + \mu^2 \left(\frac{\lambda_U^2}{2-\frac{\mu}{\delta}}\right)\\
    a_1\triangleq&\: \frac{1}{\lambda_L-\mu \frac{\lambda_U^2}{2 - \frac{\mu}{\delta}}} = O(1)\\
    a_2\triangleq&\: \frac{25 N \|\mathcal{J}_{\epsilon}\|^2 \|\mathcal{V}_{R}\|^2}{1-\|J_{\epsilon}\|} = O(1)\\
    h_1 \triangleq&\: 9(\beta^2 + \lambda_U^2) a_1 = O(1)\\
    h_2 \triangleq&\: 3 \beta^2 = O(1)\\
    h_3 \triangleq&\: 3 \nu^2 a_1 = O(1)\\
    h_4 \triangleq&\: 6 \nu^2 \lambda_U^2 a_1 = O(1)\\
    h_5 \triangleq&\: 9 \nu^2 (\lambda_U^2 + \beta^2) a_1 = O(1)\\
    h_6 \triangleq&\: 3 \nu^2 \beta^2 = O(1)\\
    h_7 \triangleq&\: 2 a_2 = O(1)\\
    h_8 \triangleq&\: \left(2\lambda_U^2 + \frac{1-\|\mathcal{J}_{\epsilon}\|}{25}3\beta^2\right) a_2 = O(1)\\
    h_9 \triangleq&\: 3 \left(\lambda_U^2 + \beta^2\right) a_2 = O(1)\\
    h_{10} \triangleq&\: \nu^2 a_2 = O(1)\\
    h_{11} \triangleq&\: \nu^2 \left(2\lambda_U^2 + \frac{1-\|\mathcal{J}_{\epsilon}\|}{25}3\beta^2\right) a_2 = O(1)\\
    h_{12} \triangleq&\: \nu^2 \left(1 + 3 \lambda_U^2 + 3 \beta^2\right) a_2 = O(1)\\
    b_1 \triangleq&\: 9 a_1 \|g(w_{\delta}^o)\|^2 = O(1)\\
    b_2 \triangleq&\: 3 a_1(3 \beta^2\|w_{\delta}^o\|^2 + \sigma^2) = O(1)\\
    b_3 \triangleq&\: 3 \beta^2\|w_{\delta}^o\|^2 + \sigma^2 = O(1)\\
    b_4 \triangleq&\: 2 a_2 \left(\delta^2 \|r(\mathds{1} \otimes w_{\delta}^o)\|^2\right) = O(1)\\
    b_5 \triangleq&\: 2 a_2 \|g(\mathds{1} \otimes w_{\delta}^o)\|^2 + \|\mathcal{J}_{\epsilon}\|^2 \|\mathcal{V}_{R}\|^2 N\left( 3 \beta^2 \|w_{\delta}^o\|^2 + \sigma^2\right) \notag \\
    =&\: O(1)\\
    b_6 \triangleq&\: a_2 \left( 3 \|g(\mathds{1} \otimes w_{\delta}^o)\|^2 + 3  \beta^2 N \|w_{\delta}^o\|^2 + N \sigma^2\right) = O(1)
  \end{align}
\end{lemma}
\begin{IEEEproof} Appendix~\ref{ap:errorrec}. \end{IEEEproof}
It is evident from expression~\eqref{eq:driving_matrix} that the stability of the driving matrix \( \Gamma \) depends critically on the fraction between the step-size \( \mu \) and the smoothing parameter \( \delta \). Motivated by this observation, let us set, for a small \( \kappa > 0 \):
\begin{equation}
  \delta = \mu^{\frac{1}{2} - \kappa}
\end{equation}
so that
\begin{equation}
  \frac{\mu}{\delta^2} = \mu^{2\kappa} \to 0 \ \mathrm{as}\ \mu \to 0
\end{equation}
Under this construction, the driving matrix satisfies
\begin{equation}
  \Gamma= \begin{bmatrix} \gamma_c + O(\mu^2) & O(\mu^{2\kappa})\\
    O(\mu^{1+2\kappa}) & \|\mathcal{J}_{\epsilon}\|+O(\mu^{1+2\kappa})
  \end{bmatrix}\\
\end{equation}
which ensures that the off-diagonal coupling terms diminish as \( \mu, \delta \to 0 \).

\begin{lemma}\label{lem:lemma8}
  Let \( \delta = \mu^{\frac{1}{2} - \kappa},\ \frac{1}{2} > \kappa > 0 \). Then there exists a small enough \( \mu \), such that \( \rho(\Gamma)< 1 \). Furthermore,
  \begin{equation}
    \underset{i \to \infty}{\limsup} \: \begin{bmatrix} \E \|\widetilde{\w}_{c,i}\|^2\\
      \E \|\bcw_{e,i}\|^2 \end{bmatrix} \preceq \begin{bmatrix} O(\mu) + O(\mu^{4\kappa}) \\ O(\mu^{1+2\kappa}) \end{bmatrix}
  \end{equation}
\end{lemma}
\begin{IEEEproof}
  See Appendix~\ref{ap:lemma8}.
\end{IEEEproof}
\noindent
\begin{theorem}\label{th:finalresult}
  Let \( \delta = \mu^{\frac{1}{2} - \kappa}, \ \frac{1}{2} > \kappa > \frac{1}{4} \). Then it holds that for sufficiently small \( \mu \),
  \begin{equation}
    \underset{i \to \infty}{\limsup} \: \E \|w_{\delta}^o - \w_{k,i}\|^2 = O(\mu)
  \end{equation}
\end{theorem}
\begin{IEEEproof}
  We have
  \begin{align}
    \E \|w_{\delta}^o - \w_{k,i}\|^2 =&\: \E \|\widetilde{\w}_{c,i} + \left(v_{L,k} \otimes I_M\right) \bcw_{e,i}\|^2\notag \\
    \le&\:2 \E \|\widetilde{\w}_{c,i}\|^2 + 2\nu^2 \E \|\bcw_{e,i}\|^2
  \end{align}
  so that the theorem follows after taking the limit and applying Lemma~\ref{lem:lemma8}.
\end{IEEEproof}

\section{Application: Division of Labor in Machine Learning}\noindent
We illustrate the performance of the algorithm in an online machine learning problem over a heterogeneous network. Given
random binary class variables \( \boldsymbol{\gamma} = \pm 1 \) and feature vectors \( \boldsymbol{h} \in \mathds{R}^M \),
the general objective in single-agent machine learning is to find a classifier \( c^{\star}(\boldsymbol{h}) \), such
that
\begin{equation}\label{eq:exactclass}
  c^{\star} \triangleq \underset{c}{\arg\min} \: \mathrm{Prob}\left \{ c(\boldsymbol{h}) \neq \boldsymbol{\gamma}
  \right \}.
\end{equation}
We restrict the class of permissible classifiers to linear classifiers of the form \( c(\boldsymbol{h}) =
\boldsymbol{h}^{\T} w \) with \( w \in \mathds{R}^M \) and approximate~\eqref{eq:exactclass} by the logistic cost to
obtain
\begin{equation}\label{eq:approximateclass}
  w^o \triangleq \underset{w}{\arg\min} \:\E \:\ln\big[ 1 + e^{-\boldsymbol{\gamma}
    \boldsymbol{h}^{\T}{w}}\big]
\end{equation}
\subsection{Group Lasso}
\noindent Regularization is an effective technique to incorporate prior structural knowledge about the classifier into the
optimization problem as a means to avoiding overfitting and improving generalization ability. For example, when the
linear classifier is known to be sparse, regularization through the \( \ell_1 \)-norm, also known as Lasso-regularization,
has been shown to encourage sparse solutions~\cite{Tibshirani94}. When there is further knowledge about the structure of the sparsity, the
group-Lasso has been proposed~\cite{Kim06, Meier08}. It takes the form
{\begin{equation}\label{eq:group}
  R(w) = \sum_{k} \lambda_k\| D_k w\|_1 = \sum_{k} \lambda_k\|w_g^k\|_1
\end{equation}
where
\begin{equation}
  w_g^k \triangleq D_k w
\end{equation}
and \( D_k \) denotes a diagonal selection matrix with entries \( 0 \) or \( 1 \) where \( 1 \)'s appear for entries of \( w \) belonging to a group. {Note that in contrast to the traditional group Lasso employing \( \ell_2 \)-norms, we are using here \( \ell_1 \)-norms to encourage within-group sparsity as well. The proposed algorithm is equally applicable to the standard group Lasso problem from~\cite{Kim06, Meier08}.}} Relation~\eqref{eq:group} is in the form of a sum-of-costs and hence
immediately decomposable.
\subsection{Network Structure}
\noindent We consider a network consisting of \( 3 \) types of agents: fully-informed (\( \mathcal{F} \)), data-informed (\(
\mathcal{D} \)), and structure-informed (\( \mathcal{S} \)) agents. Fully-informed agents have access to streaming
realizations \( \left \{ \boldsymbol{\gamma}_k(i), \boldsymbol{h}_{k,i}\right \} \) as well as knowledge about a subset
of covariates of \( w \) which are likely to be sparse, collected in \( w_{g}^k \). These agents are equipped with the
regularized cost \( J_k(w)+R_k(w) \), where
\begin{align}
  J_k(w) &= \E \:\ln\big[ 1 + e^{-\boldsymbol{\gamma}_k \boldsymbol{h}_{k}^{\T}{w}}\big] +
  \rho_2 \|w\|_2^2\\ R_k(w) &= \rho_1 \|w_{g}^k\|_1 \label{eq:fullyaug}
\end{align}
for \( k \in \mathcal{F} \). Data-informed agents have access to streaming realizations \( \left \{ \boldsymbol{\gamma}_k(i), \boldsymbol{h}_{k,i}\right \} \), but no knowledge about the structure of sparsity in \( w \). They are equipped with
\begin{align}
  J_k(w) &= \E \:\ln\big[ 1 + e^{-\boldsymbol{\gamma}_k \boldsymbol{h}_{k}^{\T}{w}}\big] + \rho_2 \|w\|_2^2\\
  R_k(w) &= 0
\end{align}
for \( k \in \mathcal{D} \). Structure-informed agents have information about the sparsity of \( w \), but no access to
realizations of feature vectors. They are equipped with
\begin{align}
  J_k(w) &= 0\\
  R_k(w) &= \rho_1 {\|w_{g}^k\|}_1\label{eq:structaug}
\end{align}
for \( k \in \mathcal{S} \).
\begin{figure}
  \includegraphics[width=.9\linewidth]{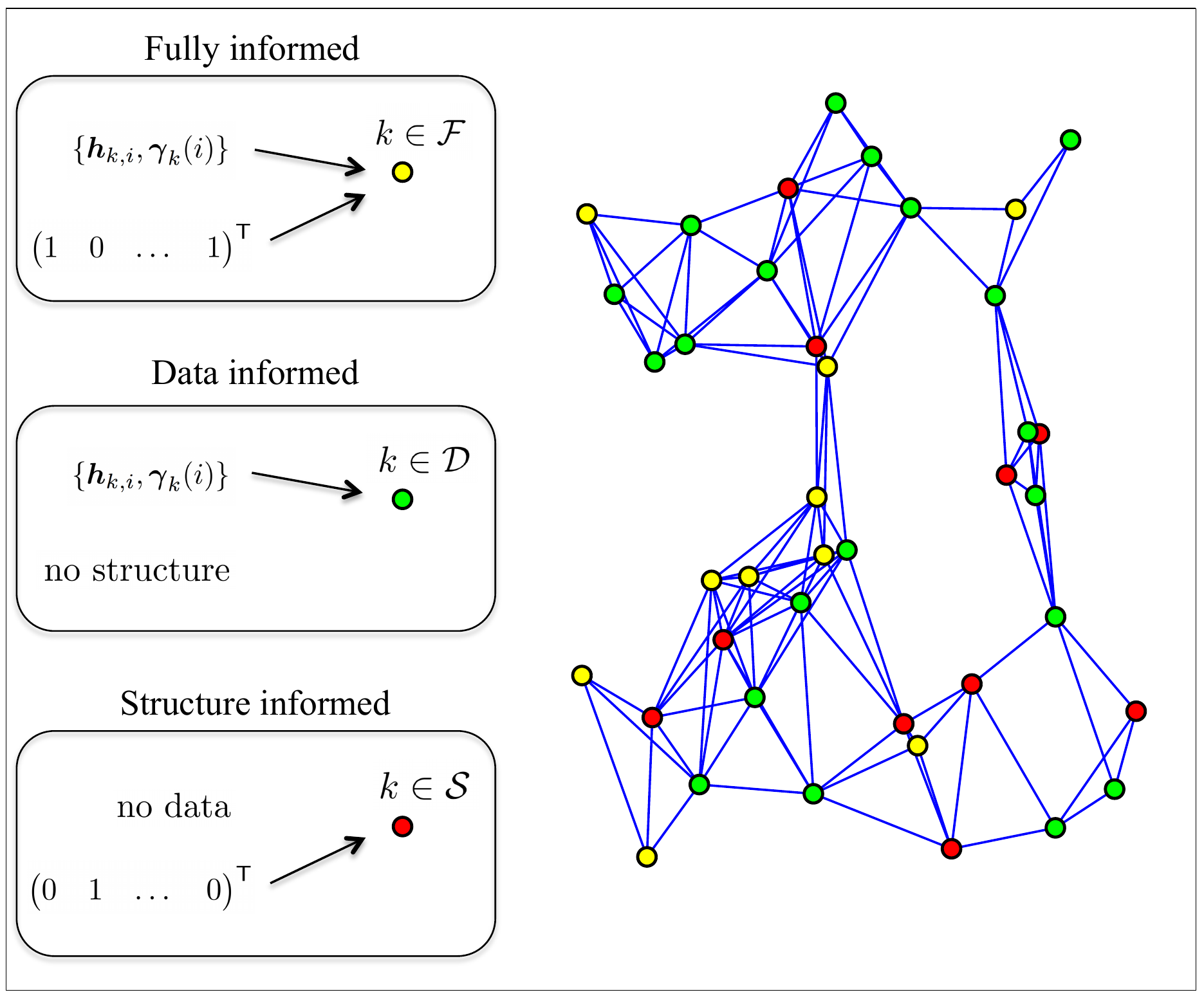}
  \centering
  \caption{Sample network consisting of \( N=40 \) agents, \( \mathrm{card}(\mathcal{F})=10 \), \( \mathrm{card}(\mathcal{D})=20 \),
\( \mathrm{card}(\mathcal{S})=10 \). Fully-informed agents have access to data as well as partial structural information.
Data-informed agents observe realizations of the feature vector along with class-labels, but have no information on the
structure of the classifier. Structure-informed agents do not have access to data, but do have partial information on
sparse elements.}\label{fig:network}
\end{figure}
Similar to ordinary \( \ell_1 \)-norm regularization, the proximal operator of \( \rho_1 \|w_g^k\|_1 \) is available in
closed form as a variation of soft-thresholding. Note that \( \|w_g^k\|_1 = \|D_k w\|_1 \), where \( D_k \) is a
diagonal matrix with \( D_{(ii)} = 1 \), if the \( i-th \) element of \( w \) is likely to be sparse and \( 0 \)
otherwise. We then obtain
\begin{equation}
  \mathrm{prox}_{\delta \rho_1 \|w_g^k\|_1}\left( w \right) = D_k \mathrm{prox}_{\delta \rho_1 \|w\|_1}\left( w \right).
\end{equation}
It is hence possible for each agent \( k \) to run~\eqref{eq:adapt1}--\eqref{eq:combine}. As long as at least one agent in
the network is either fully-informed or data-informed, the weighted sum of costs across the network is strongly convex and
assumptions~\ref{as:lipschitz_gradient_general_regularizer} through~\ref{as:reg} are satisfied. We conclude from Theorem~\ref{th:finalresult} that all agents in the network will converge to the neighborhood of:
\begin{align}
  w^o = \underset{w}{\arg\min} \: \sum_{k \in \mathcal{F}\cup\mathcal{D}} p_k \left \{ \E \: \ln\big[ 1 + e^{-\boldsymbol{\gamma}_k \boldsymbol{h}_{k}^{\T}{w}}\big]  \right \}\notag \\
  + \rho_2 \cdot \mathrm{card}(\mathcal{F}\cup\mathcal{D}) \|w\|_2^2 + \sum_{k \in \mathcal{F}\cup\mathcal{S}} p_k \|w_g^k\|_1
\end{align}
where \( \mathrm{card}(\mathcal{F}\cup\mathcal{D}) \) denotes the cardinality of the set \( \mathcal{F}\cup\mathcal{D} \), i.e. the number of agents who are either fully or data-informed. This classifier minimizes the weighted average logistic cost across the network, hence incorporating data from
all agents, regularized by the \( \ell_2 \)-norm and weighted group Lasso. Through local interactions, both data and
structural information is diffused across the entire network, allowing all agents, irrespective of their type and
available information, to arrive at an accurate classification decision.

\subsection{Numerical Results}\noindent
Performance is illustrated on the network depicted in Fig.~\ref{fig:network}, consisting of a total of \( N=40 \) agents,
\( 20 \) of which are data-informed and \( 10 \) each of which are fully and structure informed respectively. The network is
heterogeneous in both the types of available information and the noise profile of feature realizations, when data is
available. Features are generated as
\begin{equation}
  \boldsymbol{h}_{k,i} = \boldsymbol{\gamma}(i) \begin{pmatrix} 1&1&\cdots&0&0 \end{pmatrix}^{\T} + \boldsymbol{v}_k(i)
\end{equation}
where \( \boldsymbol{v}_k(i) \sim \mathcal{N}(0,\sigma_{v,k}^2) \) and \( \begin{pmatrix} 1&1&\cdots&0&0
\end{pmatrix}^{\T} \) consists of \( 50 \) leading \( 1 \)'s followed by \( 50 \) trailing \( 0 \)'s.
\begin{figure}
  \includegraphics[width=.9\linewidth]{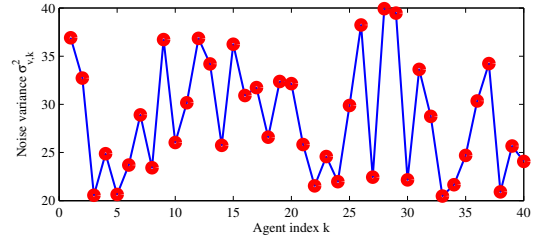}
  \centering
  \caption{Noise profile across the network for training (if \( k \in \mathcal{F}\cup\mathcal{D} \)) and testing.}\label{fig:noise}
\end{figure}
It is evident, that all class information is contained in the first half of the feature vector. This information is
dispersed across the network as follows. The noise profile across the network is depicted in Fig.~\ref{fig:noise}.

Each agent with \( k \in \mathcal{F}\cup\mathcal{S} \), i.e., fully and data-informed agents, are supplied with \( 5 \)
indices, chosen uniformly at random, of irrelevant feature covariates. They use this information to augment their cost
by an appropriate regularization as in~\eqref{eq:fullyaug} and~\eqref{eq:structaug}.
\begin{figure}
  \includegraphics[width=.9\linewidth]{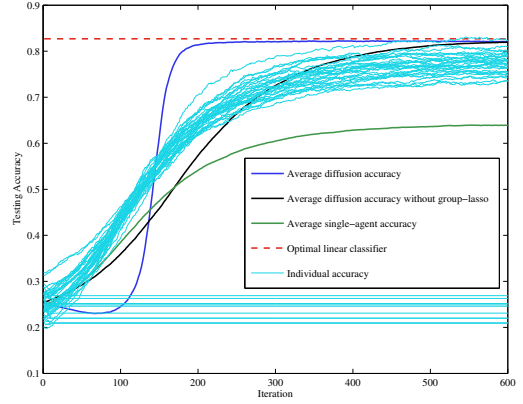}
  \centering
  \caption{Classifier performance on separate testing set.}\label{fig:performance}
\end{figure}

{The evolution of performance is illustrated in Fig.~\ref{fig:performance}. We observe that the diffusion strategy with structured sparsity regularization quickly approaches the performance of the optimal linear classifier. The rate of convergence is reduced in the absence of regularization. Finally, when no cooperation takes place, and hence information does not diffuse across the network, agents without access to observations, and those with noisy data, perform significantly worse than the cooperative strategy.}

\appendices%

\ifarx
\section{Proof of Lemma~\ref{lem:subgradients}}\label{ap:subgradients}
\noindent Let \(v \in \partial G(w)\). From the definition of the conjugate:
\begin{equation}\label{eq:inproofconj}
  G^{\star}(v) = \underset{u}{\sup} \: \left( v^{\T} u - G(u) \right)
\end{equation}
The optimality condition of the above supremum dictates that
\begin{equation}
  0 \in v - \partial G(w) \Longleftrightarrow w=\underset{u}{\arg\max} \: \left( v^{\T} u - G(u) \right)
\end{equation}
so for \(v \in \partial G(w)\), the supremum~\eqref{eq:inproofconj} is attained at \(w\). Then
\begin{equation}
  \label{eq:sumdeco}
  G^{\star}(v)=v^{\T} w - G(w).
\end{equation}
Now for any \(x\) (where the supremum might in general not be attained):
\begin{align}
  G^{\star}(x) &= \underset{u}{\sup} \: \left( x^{\T} u - G(u) \right)\notag \\
  &\ge x^{\T} w - G(w)\notag \\
  &= v^{\T} w - G(w) + w^{\T} (x-v)\notag \\
  &=G^{\star}(v) + w^{\T} (x-v)
\end{align}
By definition, any vector that satisfies \(G^{\star}(x) -G^{\star}(v) \ge w^{\T} (x-v)\) for all \(x\) is a subgradient
of \(G^{\star}(\cdot)\) at \(v\), i.e., \(w \in \partial G^{\star}(v)\). The other direction follows analogously, after
noting that for closed, convex functions \({\left(G^{\star}(\cdot)\right)}^{\star} = G(\cdot)\).
\fi
\section{Proof of Theorem~\ref{th:smoothgradient}}\label{ap:smoothgradient}
\noindent Let \(u^o \in \partial {\left(R_k^{\star} + \delta \cdot  d\right)}^{\star}(w) = \partial R_k^{\delta}(w)\). From Lemma~\ref{lem:subgradients}, this is equivalent to
\begin{equation}
  w \in \partial R_k^{\star}(u^o) + \delta \cdot  d(u^o)
\end{equation}
which due to optimality conditions is equivalent to
\begin{equation}
  u^o=\underset{u \in \dom R_k^{\star}}{\arg\max} \: \left \{ w^{\T} u - R_k^{\star}(u) - \delta\cdot d(u) \right \}.
\end{equation}
Since \(R_k^{\star}(w) + \delta\cdot d(w)\) is strongly-convex, the minimizer \(u^o\) is unique and the above holds for
any \(u^o \in \partial R_k^{\delta}(w)\). We conclude that the set \(\partial R_k^{\delta}(w)\) and hence
\begin{equation}
  \left \{ \partial R_k^{\delta}(w) \right \} = \nabla R_k^{\delta}(w) = u^o.
\end{equation}
To prove the bound on the gradient of the smooth approximation, let \(u_1^o = \nabla R_k^{\delta}(w_1)\) and
\(u_2^o=\nabla R_k^{\delta}(w_2)\) for any \(w_1, w_2\). From Lemma~\ref{lem:subgradients}, this implies \(w_1 \in
\partial R_k^{\star}(u^o_1) + \delta \cdot \partial  d(u^o_1)\) and \(w_2 \in \partial R_k^{\star}(u^o_2) + \delta
\cdot \partial  d(u^o_2)\). From the strong-convexity of \(\delta\cdot d(\cdot)\), we have:
\begin{align}
  &\:{\left(R_k^{\star}(u^o_1) + \delta \cdot \partial  d(u^o_1) - \partial R_k^{\star}(u^o_2) + \delta
  \cdot \partial  d(u^o_2)\right)}^{\T}(u^o_1 - u^o_2) \notag \\
  \ge&\: \delta \|u^o_1 - u^o_2\|^2
\end{align}
Plugging in \(w_1 \in \partial R_k^{\star}(u^o_1) + \delta \cdot \partial d(u^o_1)\) and \(w_2 \in \partial
R_k^{\star}(u^o_2) + \delta \cdot \partial  d(u^o_2)\) as well as \(u_1^o = \nabla R_k^{\delta}(w_1)\) and
\(u_2^o=\nabla R_k^{\delta}(w_2)\) yields
\begin{align}
  &\: {(w_1 - w_2)}^{\T}(\nabla R_k^{\delta}(w_1) - \nabla R_k^{\delta}(w_2)) \notag \\
  \ge&\: \delta \|\nabla R_k^{\delta}(w_1) - \nabla R_k^{\delta}(w_2)\|^2
\end{align}
which is the co-coercitivity property~\eqref{eq:rkcoco}.

\section{Proof of Theorem~\ref{th:smoothacc}}\label{ap:smoothacc}\noindent
\allowdisplaybreaks
\noindent For ease of exposition, let us introduce
\begin{align}
  F(w) \triangleq \sum_{k=1}^{N} p_k \left \{ J_k(w) + R_k(w) \right \} \label{eq:global_cost} \\
  F^{\delta}(w) \triangleq \sum_{k=1}^{N} p_k \left \{ J_k(w) + R_k^{\delta}(w) \right \} \label{eq:global_smoothed_cost}
\end{align}
We establish a string of inequalities around the difference in function values \( F(w^o) - F^{\delta}(w_{\delta}^o) \). On one hand, we have:
\begin{align}
  &\: F(w^o) - F^{\delta}(w_{\delta}^o) \notag \\
  \ifarx=&\: \sum_{k=1}^N p_k \left \{J_k(w^o) + R_k(w^o) \right \} - \sum_{k=1}^N p_k \left \{J_k(w_{\delta}^o) + R_k^{\delta}(w_{\delta}^o) \right \} \notag \\
  =&\: \sum_{k=1}^N p_k \left \{J_k(w^o) - J_k(w_{\delta}^o) \right \} + \sum_{k=1}^N p_k \left \{R_k(w^o) - R_k^{\delta}(w_{\delta}^o) \right \} \notag \\ \fi
  \stackrel{(a)}{=}&\: \sum_{k=1}^N p_k \left \{J_k(w^o) - J_k(w_{\delta}^o) \right \} \notag \\
  &\: + \sum_{k=1}^N p_k \left \{R_k(w^o) - \max_u \left( u^{\T} w_{\delta}^o - R_k^{\star}(u) - \delta d(u)\right) \right \} \notag \\
  \stackrel{(b)}{=}&\: \sum_{k=1}^N p_k \left \{J_k(w^o) - J_k(w_{\delta}^o) \right \} \notag \\
  &\: + \sum_{k=1}^N p_k \Big \{R_k(w^o) - {\nabla R_k^{\delta}(w_{\delta}^o)}^{\T} w_{\delta}^o + R_k^{\star}(\nabla R_k^{\delta}(w_{\delta}^o)) \notag \\
  &\: \ \ \ \ \ \ \ \ \ \ \ \ + \delta d\left(\nabla R_k^{\delta}(w_{\delta}^o)\right)\Big \} \notag \\
  \stackrel{(c)}{\ge}&\: \sum_{k=1}^N p_k \left \{J_k(w^o) - J_k(w_{\delta}^o) \right \} \notag \\
  +& \sum_{k=1}^N p_k \left \{\nabla {R_k^{\delta}(w_{\delta}^o)}^{\T}w^o - {\nabla R_k^{\delta}(w_{\delta}^o)}^{\T} w_{\delta}^o  + \delta d\left(\nabla R_k^{\delta}(w_{\delta}^o)\right)\right \} \notag \\
  \stackrel{(d)}{\ge}&\: \sum_{k=1}^N p_k {\nabla J_k(w_{\delta}^o)}^{\T}\left( w^o - w_{\delta}^o\right) + \frac{\lambda_L}{2}\|w^o - w_{\delta}^o\|^2 \notag \\
  +& \sum_{k=1}^N p_k \left \{{\nabla R_k^{\delta}(w_{\delta}^o)}^{\T}w^o - {\nabla R_k^{\delta}(w_{\delta}^o)}^{\T} w_{\delta}^o  + \delta d\left(\nabla R_k^{\delta}(w_{\delta}^o)\right)\right \} \notag \\
  =&\: \sum_{k=1}^N p_k {\left(\nabla J_k(w_{\delta}^o) + \nabla R_k^{\delta}(w_{\delta}^o) \right)}^{\T}\left( w^o - w_{\delta}^o\right) \notag \\
  &\: + \frac{\lambda_L}{2}\|w^o - w_{\delta}^o\|^2 + \sum_{k=1}^N p_k \delta  d\left(\nabla R_k^{\delta}(w_{\delta}^o)\right) \notag \\
  \stackrel{(e)}{=}&\: \frac{\lambda_L}{2}\|w^o - w_{\delta}^o\|^2 + \sum_{k=1}^N p_k \delta  d\left(\nabla R_k^{\delta}(w_{\delta}^o)\right)\label{eq:intermediate_1}
\end{align}
Here, \((a)\) follows from the definition of the smooth approximation~\eqref{eq:smoothshort}, \((b)\) follows from the
expression for the gradient of the smooth approximation~\eqref{eq:gradient}, \((c)\) follows from the property \(
R^{\star}(x) \triangleq \sup_u \left( u^{\T}x - R(u) \right) \ge y^{\T} x - R(y) \ \forall \ x, y\), \((d)\) follows
from the aggregate strong convexity~\eqref{eq:jkstrong} and \((e)\) follows from the definition of \(w_{\delta}^o\) and
the minimizer of the smoothed aggregate cost.

To prove the upper bound, we bound the bias for each agent individually. To begin with, note that convexity of \( J_k(\cdot) \) and \( R_k(\cdot) \) yields for all \( r_k(w^o) \in \partial R_k(w^o) \):
\begin{align}
  \ifarx J_k(w_{\delta}^{o}) - J_k(w^{o}) & \ge  {\left( \nabla J_k(w^o) \right)}^{\T} \left( w_{\delta}^o - w^o \right) \notag \\
  \Longleftrightarrow\fi J_k(w^{o}) - J_k(w_{\delta}^{o}) &\le  {\left( \nabla J_k(w^o) \right)}^{\T} \left( w^o - w_{\delta}^o \right)\\
  R_k(u) - R_k(w^{o}) &\ge  {\left( r_k(w^o) \right)}^{\T} \left( u - w^o \right)
\end{align}
Then,
\begin{align}
  &\: J_k(w^{o}) + R_k(w^o) - J_k(w_{\delta}^{o}) - R_k^{\delta}(w_{\delta}^o) \notag \\
  \ifarx=&\: J_k(w^{o}) + R_k(w^o) - J_k(w_{\delta}^{o}) \notag \\
  &\:- \min_u \left \{ R_k(u) + \delta d^{\star}\left( \frac{w_{\delta}^{o} - u}{\delta} \right) \right \} \notag \\ \fi
  =&\: J_k(w^{o}) - J_k(w_{\delta}^{o}) \notag \\
  &\: - \min_u \left \{ R_k(u) - R_k(w^o) + \delta d^{\star}\left( \frac{w_{\delta}^{o} - u}{\delta} \right) \right \} \notag \\
  \le&\: {\left( \nabla J_k(w^o) \right)}^{\T} \left( w^o - w_{\delta}^o \right) \notag \\
  &\: - \min_u \left \{ {\left( r_k(w^o) \right)}^{\T} \left( u - w^o \right) + \delta d^{\star}\left( \frac{w_{\delta}^{o} - u}{\delta} \right) \right \} \notag \\
  \ifarx=&\: {\left( \nabla J_k(w^o) +  r_k(w^o) \right)}^{\T} \left( w^o - w_{\delta}^o \right) \notag \\
  &\: - \min_u \left \{ {\left( r_k(w^o) \right)}^{\T} \left( u - w_{\delta}^o \right) + \delta d^{\star}\left( \frac{w_{\delta}^{o} - u}{\delta} \right) \right \} \notag \\ \fi
  \stackrel{(a)}{=}&\: {\left( \nabla J_k(w^o) +  r_k(w^o) \right)}^{\T} \left( w^o - w_{\delta}^o \right) \notag \\
  &\: - \delta \min_v \left \{ - {\left( r_k(w^o) \right)}^{\T} v + d^{\star}\left( v \right) \right \} \notag \\
  \ifarx=&\: {\left( \nabla J_k(w^o) +  r_k(w^o) \right)}^{\T} \left( w^o - w_{\delta}^o \right) \notag \\
  &\: + \delta \max_v \left \{ {\left( r_k(w^o) \right)}^{\T} v - d^{\star}\left( v \right) \right \} \notag \\ \fi
  \stackrel{(b)}{=}&\: {\left( \nabla J_k(w^o) +  r_k(w^o) \right)}^{\T} \left( w^o - w_{\delta}^o \right) + \delta d\left( r_k(w^o) \right)
\end{align}
where \( (a) \) follows after a change of variables \( v \triangleq \frac{w_{\delta}^{o} - u}{\delta} \) and \( (b) \) is a result of the definition of the conjugate function. Returning to the aggregate cost, we then have
\begin{align}
  &\: \sum_{k=1}^{N} p_k \left \{ J_k(w^{o}) + R_k(w^o) \right \} -\sum_{k=1}^{N} p_k \left \{ J_k(w_{\delta}^{o}) + R_k^{\delta}(w_{\delta}^o) \right \} \notag \\
  \ifarx =&\: \sum_{k=1}^{N} p_k \left \{ J_k(w^{o}) + R_k(w^o) - J_k(w_{\delta}^{o}) + R_k^{\delta}(w_{\delta}^o) \right \} \notag \\ \fi
  \ifarx\le&\: \sum_{k=1}^{N} p_k \left \{ {\left( \nabla J_k(w^o) +  r_k(w^o) \right)}^{\T} \left( w^o - w_{\delta}^o \right) \right \} \notag \\\fi
  &\: + \sum_{k=1}^{N} p_k \delta d\left( r_k(w^o) \right) \notag \\
  \ifarx = \else \le \fi&\: {\left \{ \sum_{k=1}^{N} p_k \left( \nabla J_k(w^o) +  r_k(w^o) \right) \right \}}^{\T} \left( w^o - w_{\delta}^o \right) \notag \\
  &\: + \sum_{k=1}^{N} p_k \delta d\left( r_k(w^o) \right)
\end{align}
By definition, \( w^o \) is the minimizer of \( \sum_{k=1}^{N} p_k \left \{ J_k(w^{o}) + R_k(w^o) \right \} \), so there exist subgradients \( r_k^o \in \partial R_k(w^o) \), such that
\begin{equation}
  \sum_{k=1}^N p_k \left( \nabla J_k(w^o) +  r_k^o \right) = 0
\end{equation}
Then,
\begin{align}
  &\: \sum_{k=1}^{N} p_k \left \{ J_k(w^{o}) + R_k(w^o) \right \} -\sum_{k=1}^{N} p_k \left \{ J_k(w_{\delta}^{o}) + R_k^{\delta}(w_{\delta}^o) \right \} \notag \\
  \le&\: \sum_{k=1}^{N} p_k \delta d\left( r_k^o \right) = O(\delta)
\end{align}
\ifarx We conclude from~\eqref{eq:intermediate_1}:
\begin{align}
  &\: \frac{\lambda_L}{2}\|w^o - w_{\delta}^o\|^2 + \sum_{k=1}^N p_k \delta  d\left(\nabla R_k^{\delta}(w_{\delta}^o)\right) \notag \\
  \le&\:  F(w^o) - F^{\delta}(w_{\delta}^o) \le  \sum_{k=1}^{N} p_k \delta d\left( r_k^o \right)
\end{align}
The result follows after rearranging. \else
The result follows after combining with~\eqref{eq:intermediate_1} and rearranging. \fi

\section{Proof of Lemma~\ref{lem:centcontr}}\label{ap:centcontr}
\noindent Let \(\alpha \) be an arbitrary real number such that \(0<\alpha<1\). Then
\begin{align}
  \phantom{=}& \|T_c(x)-T_c(y)\|^2 \notag \\
  =& \Big \| x - y - \mu \sum_{k=1}^{N} p_k \Big \{ \nabla J_k(x) - \nabla J_k(y) \notag \\
  &\: \ \ \ \  \ \ \ \ \ \ \ \ \ \ \ \ \ \ \ \ \ \ \ + \nabla R_k^{\delta}(x) - \nabla R_k^{\delta}(y) \Big \} \Big \|^2\notag \\
  =& \|x-y\|^2 + \mu^2 \Bigg \| \sum_{k=1}^{N} p_k \Big \{ \nabla J_k(x) - \nabla J_k(y) \notag \\
  &\: \ \ \ \ \ \ \ \ \ \ \ \ \ \ \ \ \ \ \ \ \ \ \ \ \ \ \ \ \  + \nabla R_k^{\delta}(x) - \nabla R_k^{\delta}(y) \Big \} \Bigg \|^2\notag \\
  &- 2 \mu\sum_{k=1}^{N} p_k {\left(x-y\right)}^{\T} \left(\nabla J_k(x) - \nabla J_k(y)\right) \notag \\
  &\: - 2 \mu\sum_{k=1}^{N} p_k {\left(x-y\right)}^{\T} \left(\nabla R_k^{\delta}(x) - \nabla R_k^{\delta}(y)\right)\notag \\
  \stackrel{(a)}{\le}& \|x-y\|^2 + \mu^2 \sum_{k=1}^{N} p_k \Big \|\nabla J_k(x) - \nabla J_k(y) \notag \\
  &\: \ \ \ \ \ \ \ \ \ \ \ \ \ \ \ \  \ \ \ \ \ \ \ \ \ \ + \nabla R_k^{\delta}(x) - \nabla R_k^{\delta}(y) \Big \|^2 \notag \\
  &- 2 \mu \lambda_L {\left \|x-y\right \|}^2 - 2 \mu \delta \sum_{k=1}^{N} p_k\left \| \nabla R_k^{\delta}(x) - \nabla R_k^{\delta}(y)\right \|^2\notag \\
  \stackrel{(b)}{\le}& \|x-y\|^2 + \mu^2 \sum_{k=1}^{N} p_k \frac{1}{\alpha}\Big \|\nabla J_k(x) - \nabla J_k(y)\Big \|^2 \notag \\
  &+\mu^2 \sum_{k=1}^{N} p_k \frac{1}{1-\alpha}\Big \|\nabla R_k^{\delta}(x) -  \nabla R_k^{\delta}(y) \Big \|^2 - 2 \mu \lambda_L \left \|x-y\right \|^2  \notag \\
  & - 2 \mu \delta \sum_{k=1}^{N} p_k \left \|\nabla R_k^{\delta}(x) -  \nabla R_k^{\delta}(y)\right \|^2
\end{align}
where \((a)\) follows from Jensen's inequality, strong convexity~\eqref{eq:jkstrong}, and co-coercivity~\eqref{eq:rkcoco},
and \((b)\) from \( \|a+b\|^2 \le \frac{1}{\alpha}\|a\|^2 + \frac{1}{1-\alpha}\|b\|^2 \) for any \( a,b \in \mathds{R}^M \). Since, by assumption, \( \mu<2\delta \), we select \( \alpha = 1-\frac{\mu}{2\delta} \).
This results in \( \frac{\mu^2}{1-\alpha} = 2\mu \delta \) and allows us to cancel all terms involving \( \nabla_w
R_k^{\delta}(\cdot) \) in the above inequality. Hence,
\begin{align}
  &\: \|T_c(x)-T_c(y)\|^2 \notag \\
  \le&\: \|x-y\|^2 + \mu^2 \sum_{k=1}^{N} p_k \frac{1}{1-\frac{\mu}{2\delta}}\Big \|\nabla J_k(x) - \nabla J_k(y)\Big \|^2 \notag \\
  &\: - 2 \mu \lambda_L \left \|x-y\right \|^2 \notag \\
  \stackrel{(a)}{\le}&\: \|x-y\|^2 + \frac{\mu^2 \lambda_U^2}{1-\frac{\mu}{2\delta}}\|x-y\|^2 - 2 \mu \lambda_L \left \|x-y\right \|^2\notag \\
  =&\: \left( 1-2 \mu \lambda_L + \mu^2 \frac{\lambda_U^2}{1-\frac{\mu}{2\delta}} \right)\|x-y\|^2\notag \\
  \stackrel{(b)}{\le}&\: {\left( 1- \mu \lambda_L + \mu^2 \frac{\lambda_U^2}{2-\frac{\mu}{\delta}} \right)}^2\|x-y\|^2
\end{align}
where \((a)\) is due to the Lipschitz property~\eqref{eq:gradbound} and \((b)\) is due to \(1-a \le
{(1-\frac{1}{2}a)}^2\) for all \(a \in \mathds{R}\). From Banach's fixed-point theorem, we know that as long as
\(\gamma_c<1\), \(w_i=T_c(w_{i-1})\) converges exponentially to a unique fixed point, which satisfies
\(w_{\infty}=T_c(w_{\infty})\). From~\eqref{eq:Tc}, we conclude that
\begin{equation}
  \sum_{k=1}^{N} p_k \nabla J_k(w_{\infty}) + \sum_{k=1}^{N} p_k \nabla R_k^{\delta}(w_{\infty}) = 0
\end{equation}
so that from~\eqref{eq:smoothproblem}, \(w_{\infty} = w_{\delta}^o\).

\section{Proof of Lemma~\ref{lem:pertbounds}}\label{ap:pertbounds}
\allowdisplaybreaks%
\noindent The proof of the first three inequalities relies on the Lipschitz properties of the gradients and the decomposition~\eqref{eq:centrecursion}--\eqref{eq:errorrecursion}. First, we bound the terms arising from the disagreement across the network. Denote the \(k\)-th element of
\(P[\cdot]\) by \(P_{(k)}[\cdot]\). Then
\begin{align}
  &P_{(k)}[\boldsymbol{t}_{i-1}] \notag \\
  =&\: \E \|\nabla J_k(\w_{c,i-1}) - \nabla J_k(\w_{k,i-1})\notag \\
  &\ \ \ \ \ \ + \nabla R_k^{\delta}(\w_{c,i-1}-\mu\nabla J_k(\w_{c,i-1})) \notag \\
  &\ \ \ \ \ \ - \nabla R_k^{\delta}(\w_{k,i-1}-\mu\nabla J_k(\w_{k,i-1}))\|^2\notag \\
  \stackrel{(a)}{\le}&\: 2 \E \|\nabla J_k(\w_{c,i-1}) - \nabla J_k(\w_{k,i-1})\|^2\notag \\
  &+ 2 \E \|\nabla R_k^{\delta}(\w_{c,i-1}-\mu\nabla J_k(\w_{c,i-1})) \notag \\
  &\:  \ \ \ \ \ \ \ \ \ - \nabla R_k^{\delta}(\w_{k,i-1}-\mu\nabla J_k(\w_{k,i-1}))\|^2\notag \\
  \stackrel{(b)}{\le}&\: 2 \lambda_U^2\E \|\w_{c,i-1} - \w_{k,i-1}\|^2 \notag \\
  &\: + \frac{2}{\delta^2} \E \| \w_{c,i-1}-\mu\nabla J_k(\w_{c,i-1}) \notag \\
  &\: \ \ \ \ \ \ \ \ \ \ \   - \w_{k,i-1}+\mu\nabla J_k(\w_{k,i-1})\|^2 \notag \\
  \stackrel{(c)}{\le}&\: \left(2 \lambda_U^2 + 4\frac{1 + \mu^2}{\delta^2} \right)\E \|\w_{c,i-1} - \w_{k,i-1}\|^2\notag \\
  \stackrel{(d)}{=}&\: \left(2 \lambda_U^2 + 4\frac{1 + \mu^2}{\delta^2} \right) \E \|\left(v_{L,k} \otimes I_M\right)\bcw_{e,i-1}\|^2\notag \\
  \ifarx \le&\: \left(2 \lambda_U^2 + 4\frac{1 + \mu^2}{\delta^2} \right) \nu^2 \E \|\bcw_{e,i-1}\|^2\notag \\ \fi
  =&\: \left(2 \lambda_U^2 + 4\frac{1 + \mu^2}{\delta^2} \right) \nu^2 \mathds{1}^{\T} P[\bcw_{e,i-1}]
\end{align}
where \((a)\) is due Jensen's inequality, \((b)\) and \((c)\) are due to Lipschitz continuity of the gradients and
\((d)\) is due to \(\bcw_i = \mathds{1} \otimes \w_{c,i} + \mathcal{V}_L \bcw_{e,i}\). Stacking both sides of the above
inequality yields~\eqref{eq:boundt}.

Now consider \(\boldsymbol{u}_{i-1}\), which arises from the incremental implementation:
\begin{align}
  &\:P_{(k)}[\boldsymbol{u}_{i-1}]\notag \\
  =&\: \E \|\nabla R^{\delta}_k(\w_{k,i-1}) - \nabla R_k^{\delta}(\w_{k,i-1} - \mu \nabla J_k(\w_{k,i-1}))\|^2\notag \\
  \stackrel{(a)}{\le}&\:\frac{\mu^2}{\delta^2} \E \|\nabla J_k(\w_{k,i-1})\|^2\notag \\
  =&\: \frac{\mu^2}{\delta^2} \E \|\nabla J_k(\w_{k,i-1}) - \nabla J_k(\w_{c,i-1})+ \nabla J_k( \w_{c,i-1}) \notag \\
  &\: \ \ \ \ \ \ \ \  - \nabla J_k(w_{\delta}^o) + \nabla J_k(w_{\delta}^o)\|^2\notag \\
  \stackrel{(b)}{\le}&\: \frac{\mu^2}{\delta^2} \big(3\lambda_U^2 \nu^2 \mathds{1}^{\T}P[\bcw_{e,i-1}] + 3\lambda_U^2 \E \|\widetilde{\w}_{c,i-1}\|^2 \notag \\
  &\: \ \ \ \ \  + 3 \|\nabla J_k(w_{\delta}^o)\|^2\big)
\end{align}
where \((a)\) is due to Lipschitz continuity of \(\nabla R_k^{\delta}(w)\) and \((b)\) is due to Jensen's inequality and
Lipschitz continuity of \(\nabla J_k(w)\). Upon stacking we obtain~\eqref{eq:boundu}.

Next, we bound the perturbations caused by the gradient noise \(\s_{k,i}(\w_{k,i}) = \widehat{\nabla
J}_k(\w_{k,i-1}) - \nabla  J_k(\w_{k,i-1})\). While a loose upper bound can be obtained immediately from Jensen's
inequality, it turns out that the incremental implementation~\eqref{eq:adapt2} along with the
co-coercivity~\eqref{eq:rkcoco} of \(\nabla R_k^{\delta}(w)\) have a variance reducing effect on the recursion:
\begin{align}
  &P_{(k)}[\s^g_{i} + \s^p_{i} - \E \s^p_i]\notag \\
  \stackrel{(a)}{\le}&P_{(k)}[\s^g_{i} + \s^p_{i}]\notag \\
  =&\E \left \|\nabla J_k(\w_{k,i-1}) - \widehat{\nabla J}_k(\w_{k,i-1})\right \|^2 \notag \\
  &+\E \Big \|\nabla R_k^{\delta}(\w_{k,i-1}-\mu \nabla J_k(\w_{k,i-1})) \notag \\
  &\: \ \ \ \ \ \ \ -\nabla R_k^{\delta}(\w_{k,i-1}-\mu\widehat{\nabla J}_k(\w_{k,i-1}))\Big \|^2\notag \\
  &+2\E {\Big(\nabla J_k(\w_{k,i-1}) - \widehat{\nabla J}_k(\w_{k,i-1})\Big)}^{\T}\notag \\
  &\ \ \ \ \ \ \ \times\Big(\nabla R_k^{\delta}(\w_{k,i-1}-\mu \nabla J_k(\w_{k,i-1})) \notag \\
  &\ \ \ \ \ \ \ \ \ \ \ \ \  - \nabla R_k^{\delta}(\w_{k,i-1}-\mu\widehat{\nabla J}_k(\w_{k,i-1}))\Big)\notag \\
  =&\E \left \|\nabla J_k(\w_{k,i-1}) - \widehat{\nabla J}_k(\w_{k,i-1})\right \|^2 \notag \\
  &+\E \Big \|R_k^{\delta}(\w_{k,i-1}-\mu \nabla J_k(\w_{k,i-1})) \notag \\
  &\: \ \ \ \ \ \ \ -\nabla R_k^{\delta}(\w_{k,i-1}-\mu\widehat{\nabla J}_k(\w_{k,i-1}))\Big \|^2\notag \\
  &-\frac{2}{\mu}\E \Big(\w_{k,i-1}-\mu \nabla J_k(\w_{k,i-1}) \notag \\
  &\: \ \ \ \ \ \ \ \ \ \ \ -{\big(\w_{k,i-1}-\mu \widehat{\nabla J}_k(\w_{k,i-1})\big)}\Big)^{\T}\notag \\
  &\ \ \ \times\Big(\nabla R_k^{\delta}(\w_{k,i-1}-\mu \nabla J_k(\w_{k,i-1})) \notag \\
  &\: \ \ \ \ \ \ \ \ -\nabla R_k^{\delta}(\w_{k,i-1}-\mu\widehat{\nabla J}_k(\w_{k,i-1}))\Big)\notag \\
  \stackrel{(b)}{\le}&\E \left \|\nabla J_k(\w_{k,i-1}) - \widehat{\nabla J}_k(\w_{k,i-1})\right \|^2 \notag \\
  &+\E \Big \|R_k^{\delta}(\w_{k,i-1}-\mu \nabla J_k(\w_{k,i-1})) \notag \\
  &\: \ \ \ \ \ \ \ -\nabla R_k^{\delta}(\w_{k,i-1}-\mu\widehat{\nabla J}_k(\w_{k,i-1}))\Big \|^2\notag \\
  &-\frac{2\delta}{\mu}\E \Big \|\nabla R_k^{\delta}(\w_{k,i-1}-\mu \nabla J_k(\w_{k,i-1})) \notag \\
  &\: \ \ \ \ \ \ \ \ \ \ \ -\nabla R_k^{\delta}(\w_{k,i-1}-\mu\widehat{\nabla J}_k(\w_{k,i-1}))\Big \|^2\notag \\
  \ifarx=&\E \left \|\nabla J_k(\w_{k,i-1}) - \widehat{\nabla J}_k(\w_{k,i-1})\right \|^2 \notag \\
  &- \left(\frac{2\delta}{\mu}-1\right)\E \Big \|R_k^{\delta}(\w_{k,i-1}-\mu \nabla J_k(\w_{k,i-1})) \notag \\
  &\: \ \ \ \ \ \ \ \ \ \ \ \ \ \ \ \ \ \ \ \ - \nabla R_k^{\delta}(\w_{k,i-1}-\mu\widehat{\nabla J}_k(\w_{k,i-1}))\Big \|^2\notag \\ \fi
  \stackrel{(c)}{\le}&\:\E \left \|\nabla J_k(\w_{k,i-1}) - \widehat{\nabla J}_k(\w_{k,i-1})\right \|^2\notag \\
  \ifarx =&\E \|\s_{k,i}(\w_{k,i-1})\|^2\notag \\ \fi
  \stackrel{(d)}{\le}&\beta^2\E \|\w_{k,i-1}\|^2 + \sigma^2
\end{align}
where \((a)\) follows from \(\E \|\boldsymbol{x}-\E \boldsymbol{x}\|^2 \le \E \|\boldsymbol{x}\|^2\) for any
\(\boldsymbol{x}\), \((b)\) follows from co-coercitivity~\eqref{eq:rkcoco}, \((c)\) follows from \(\mu < 2\delta \) and
\((d)\) is due to~\eqref{eq:noisebound}. Now from \(\w_{k,i-1} = \w_{c,i-1} + \left(v_{L,k} \otimes I\right)
\bcw_{e,i-1}\), we can
bound
\begin{align}
  &\: \|\w_{k,i-1}\|^2 \notag \\
  =&\: \|\w_{c,i-1} + \left(v_{L,k} \otimes I\right) \bcw_{e,i-1}\|^2\notag \\
  =&\: \|\w_{c,i-1} - w_{\delta}^o + \left(v_{L,k} \otimes I\right) \bcw_{e,i-1} + w_{\delta}^o\|^2\notag \\
  \le&\: 3\|\widetilde{\w}_{c,i-1}\|^2 + 3\nu^2 \mathds{1}^{\T}P[\bcw_{e,i-1}] + 3 \|w_{\delta}^o\|^2.
\end{align}
where we appealed to Jensen's inequality again. Eq.~\eqref{eq:bounds} follows after stacking. Next, note that because
\( \| \E \boldsymbol{x}\|^2 \le \E \|\boldsymbol{x}\|^2\)
\begin{equation}
  P[\E \s^p_i] \preceq P[\s^p_i].
\end{equation}
Subsequently,
\begin{align}
  P_{(k)}[\s^p_i] &= \E \|\nabla R_k^{\delta}(\w_{k,i-1}-\mu \nabla J_k(\w_{k,i-1})) \notag \\
  &\: \ \ \ \ \ \ \ - \nabla R_k^{\delta}(\w_{k,i-1}-\mu\widehat{\nabla J}_k(\w_{k,i-1}))\|^2\notag \\
  \stackrel{(a)}{\le}&\:\frac{\mu^2}{\delta^2}\E \|\nabla J_k(\w_{k,i-1}) - \widehat{\nabla J}_k(\w_{k,i-1})\|^2\notag \\
  =&\:\frac{\mu^2}{\delta^2}\|\s_{k,i}(\w_{k,i-1})\|^2
\end{align}
where \((a)\) is due to~\eqref{eq:rklip}, so that similarly to the above
\begin{align}
  P[\E \s^p_i] \preceq&\: 3\beta^2 \frac{\mu^2}{\delta^2} P[\mathds{1}\otimes \widetilde{\w}_{c,i-1}]+ 3 \beta^2 \frac{\mu^2}{\delta^2} \nu^2 \mathds{1}\mathds{1}^{\T}P[\bcw_{e,i-1}] \notag \\
  &\:+ 3\beta^2 \frac{\mu^2}{\delta^2}P[\mathds{1}\otimes w_{\delta}^o] + \frac{\mu^2}{\delta^2}\sigma^2 \mathds{1}
\end{align}
which is~\eqref{eq:bounde}. Next,
\begin{align}
  &\: P_{(k)}[{g}(\mathds{1} \otimes \w_{c,i-1})] \notag \\
  \ifarx =&\: \E \|\nabla J_k(\w_{c,i-1})\|^2\notag \\ \fi
  =&\: \E \|\nabla J_k(\w_{c,i-1}) - \nabla J_k(w_{\delta}^o) + \nabla J_k(w_{\delta}^o)\|^2\notag \\
  \le&\: 2 \lambda_U^2 \E \|\w_{c,i-1} - w_{\delta}^o\|^2 + 2 \|\nabla J_k(w_{\delta}^o)\|^2
\end{align}
which implies~\eqref{eq:boundgc} after stacking. Eq.~\eqref{eq:boundrc} follows analogously.

\section{Proof of Lemma~\ref{lem:errorrec}}\label{ap:errorrec}
\noindent We make use of Jensen's inequality \( \|x+y\|^2\le\frac{1}{\alpha}\|x\|^2+\frac{1}{1-\alpha}\|y\|^2\) for all
\(x,y\) and \(0<\alpha<1\):
\begin{align}
  \phantom{=}&\: \E \|\widetilde{\w}_{c,i}\|^2\notag \\
  =&\: \E \Big \|T_c(\w_{c,i-1}) - T_c(w_{\delta}^o) \notag \\
  &\: \ \ \ \ \ + \mu \left( p^{\T} \otimes I_M\right)\big( \boldsymbol{t}_{i-1} + \boldsymbol{u}_{i-1} + \s_i - \E \s_i + \E \s_i \big)\Big \|^2\notag \\
  \stackrel{(a)}{=}&\:\E \Big \|T_c(\w_{c,i-1})- T_c(w_{\delta}^o) \notag \\
  &\: \ \ \ \ \ + \mu \left( p^{\T} \otimes I_M\right)\big( \boldsymbol{t}_{i-1} + \boldsymbol{u}_{i-1} + \E \s_i \big)\Big \|^2 \notag \\
  &\:+ \mu^2 \E \Big \|\left( p^{\T} \otimes I_M\right)\big( \s_i - \E \s_i\big)\Big \|^2\notag \\
  \stackrel{(b)}{\le}&\: \frac{1}{\gamma_c}\E \Big \|T_c(\w_{c,i-1})- T_c(w_{\delta}^o)\Big \|^2 \notag \\
  &\: +\frac{\mu^2}{1-\gamma_c}\E \Big \| \left( p^{\T} \otimes I_M\right)\big( \boldsymbol{t}_{i-1} + \boldsymbol{u}_{i-1} + \E \s_i \big)\Big \|^2 \notag \\
  &\:+ \mu^2 \E \Big \|\left( p^{\T} \otimes I_M\right)\big( \s_i - \E \s_i\big)\Big \|^2\notag \\
  \stackrel{(c)}{\le}&\: {\gamma_c}\E \left \|\widetilde{\w}_{c,i-1}\right \|^2 \notag \\
  &\: + \frac{\mu^2}{1-\gamma_c}\E \Big \| \left( p^{\T} \otimes I_M\right)\big( \boldsymbol{t}_{i-1} + \boldsymbol{u}_{i-1} + \E \s_i \big)\Big \|^2 \notag \\
  &\:+ \mu^2 \E \Big \|\left( p^{\T} \otimes I_M\right)\big( \s_i - \E \s_i\big)\Big \|^2\notag \\
  \stackrel{(d)}{\le}&\: {\gamma_c}\E \left \|\widetilde{\w}_{c,i-1}\right \|^2 + \frac{\mu^2}{1-\gamma_c} p^{\T} P[ \boldsymbol{t}_{i-1} + \boldsymbol{u}_{i-1} + \E \s_i ] \notag \\
  &\: + \mu^2 p^{\T} P[\s_i - \E \s_i]\notag \\
  \stackrel{(e)}{\le}&\: {\gamma_c}\E \left \|\widetilde{\w}_{c,i-1}\right \|^2 + \frac{3\mu^2}{1-\gamma_c} p^{\T}\Big( P[ \boldsymbol{t}_{i-1}] + P[\boldsymbol{u}_{i-1}] + P[\E \s_i ]\Big) \notag \\
  &\: + \mu^2 p^{\T} P[\s_i - \E \s_i]\notag \\
  \stackrel{(f)}{\le}&\: {\gamma_c}\E \left \|\widetilde{\w}_{c,i-1}\right \|^2 + \frac{3\mu^2}{1-\gamma_c} p^{\T}\Big( \Big( 2\lambda_U^2 + \frac{1+\mu^2}{\delta^2}\Big) \nu^2 \mathds{1}\mathds{1}^{\T} P[\bcw_{e,i-1}]\notag \\
  &\: \ \ \ + \frac{\mu^2}{\delta^2} \big(3\lambda_U^2 \nu^2 \mathds{1}\mathds{1}^{\T}P[\bcw_{e,i-1}] + 3\lambda_U^2 P[\mathds{1}\otimes \widetilde{\w}_{c,i-1}] \notag \\
  &\: \ \ \ \ \ \ \ \ \ \ \  + 3 P[g(\mathds{1} \otimes w_{\delta}^o)]\big) \notag \\
  &\: \ \ \ + 3\beta^2 \frac{\mu^2}{\delta^2} P[\mathds{1}\otimes \widetilde{\w}_{c,i-1}] + 3 \beta^2 \frac{\mu^2}{\delta^2} \nu^2 \mathds{1}\mathds{1}^{\T}P[\bcw_{e,i-1}] \notag \\
  &\: \ \ \ + 3\beta^2 \frac{\mu^2}{\delta^2}P[\mathds{1}\otimes w_{\delta}^o] + \frac{\mu^2}{\delta^2}\sigma^2 \mathds{1}\Big) \notag \\
  &\:+ \mu^2 p^{\T} \Big(3\beta^2 P[\mathds{1}\otimes \widetilde{\w}_{c,i-1}] + 3 \beta^2 \nu^2 \mathds{1}\mathds{1}^{\T}P[\bcw_{e,i-1}] \notag \\
  &\: \ \ \ \ \ \ \ \ \ \ \ + 3\beta^2 P[\mathds{1}\otimes w_{\delta}^o] + \sigma^2 \mathds{1}\Big)\notag \\
  \stackrel{(g)}{=}&\: {\gamma_c}\E \left \|\widetilde{\w}_{c,i-1}\right \|^2 + \frac{3\mu^2}{1-\gamma_c} \Big( \Big(2 \lambda_U^2 + \frac{1+\mu^2}{\delta^2}\Big) \nu^2 \E \|\bcw_{e,i-1}\|^2 \notag \\
  &\: \ \ \ + \frac{\mu^2}{\delta^2} \big(3\lambda_U^2 \nu^2 \E \|\bcw_{e,i-1}\|^2 + 3\lambda_U^2 \E \|\widetilde{\w}_{c,i-1}\|^2 \notag \\
  &\: \ \ \ \ \ \ \ \ \ \ \ + 3 \|g(w_{\delta}^o)\|^2\big) \notag \\
  &\: \ \ \ + 3\beta^2 \frac{\mu^2}{\delta^2} \E \|\widetilde{\w}_{c,i-1}\|^2 + 3 \beta^2 \frac{\mu^2}{\delta^2} \nu^2 \E \|\bcw_{e,i-1}\|^2 \notag \\
  &\: \ \ \ + 3\beta^2 \frac{\mu^2}{\delta^2}\|w_{\delta}^o\|^2 + \frac{\mu^2}{\delta^2}\sigma^2\Big) \notag \\
  &\:+ \mu^2 \Big(3\beta^2 \E \|\widetilde{\w}_{c,i-1}\|^2 + 3 \beta^2 \nu^2 \E \|\bcw_{e,i-1}\|^2 \notag \\
  &\: \ \ \ \ \ \ \ \ + 3\beta^2\|w_{\delta}^o\|^2 +\sigma^2\Big)\notag \\
  \stackrel{(h)}{=}&\: \left(\gamma_c + \frac{9\mu^4(\beta^2 + \lambda_U^2)}{(1-\gamma_c)\delta^2} + 3\mu^2\beta^2\right)\E \left \|\widetilde{\w}_{c,i-1}\right \|^2 \notag \\
  &\:+\Bigg( \frac{3\mu^2\nu^2}{1-\gamma_c}\left( 2 \lambda_U^2 + \frac{1 + \mu^2 + 3 \mu^2 \lambda_U^2 + 3 \mu^2 \beta^2}{\delta^2}\right) \notag \\
  &\: \ \ \ \ \ \ +3\mu^2\beta^2 \nu^2 \Bigg) \E \|\bcw_{e,i-1}\|^2 \notag \\
  &\:+ \frac{3\mu^4}{(1-\gamma_c)\delta^2} \left(3 \|g(w_{\delta}^o)\|^2+ 3\beta^2 \|w_{\delta}^o\|^2 + \sigma^2\right) \notag \\
  &\: + \mu^2 \left(3\beta^2\|w_{\delta}^o\|^2 +\sigma^2\right)\notag \\
  \stackrel{(i)}{=}&\: \left(\gamma_c + \frac{\mu^3}{\delta^2}\frac{9(\beta^2 + \lambda_U^2)}{\lambda_L-\mu \frac{\lambda_U^2}{2-\frac{\mu}{\delta}}} + 3\mu^2 \beta^2\right)\E \left \|\widetilde{\w}_{c,i-1}\right \|^2 \notag \\
  &\:+\Bigg(\frac{\mu}{\delta^2} \frac{3\nu^2}{\lambda_L - \mu \frac{\lambda_U^2}{2-\frac{\mu}{\delta}}} + \mu \frac{6\nu^2}{\lambda_L - \mu \frac{\lambda_U^2}{2-\frac{\mu}{\delta}}} \lambda_U^2 \notag \\
  &\: \ \ \ + \frac{\mu^3}{\delta^2}\frac{9\nu^2}{\lambda_L - \mu \frac{\lambda_U^2}{2-\frac{\mu}{\delta}}}(\lambda_U^2+\beta^2) + 3\mu^2 \beta^2 \nu^2 \Bigg) \E \|\bcw_{e,i-1}\|^2 \notag \\
  &\:+ \frac{\mu^3}{\delta^2}\frac{9}{\lambda_L - \mu \frac{\lambda_U^2}{2-\frac{\mu}{\delta}}} \|g(w_{\delta}^o)\|^2 \notag \\
  &\: + \frac{\mu^3}{\delta^2}\frac{3}{\lambda_L - \mu \frac{\lambda_U^2}{2-\frac{\mu}{\delta}}} \left(3 \beta^2 \|w_{\delta}^o\|^2 + \sigma^2\right) \notag \\
  &\: + \mu^2 \left(3\beta^2\|w_{\delta}^o\|^2 +\sigma^2\right)
\end{align}
In step (a), cross-terms are eliminated because \( \E \left \{ \s_i - \E \s_i\right \}=0\). Step (b) is due to
\(\gamma_c < 1\) and Jensen's inequality, (c) is due to Lemma~\ref{lem:centcontr}, (d) and (e) follow from Jensen's
inequality. The bounds from Lemma~\ref{lem:pertbounds} are used in (f) and (g) is due to \( \mathds{1}^{\T}
P[\boldsymbol{x}]=\E \|\boldsymbol{x}\|^2\) for \(\boldsymbol{x} \in \mathds{R}^{MN}\) and \(p^{\T}
P[\mathds{1}\otimes\boldsymbol{y}]=\E \|\boldsymbol{y}\|^2\) for \( \boldsymbol{y} \in \mathds{R}^{M}\). In \((i)\), the terms are rearranged to expose the dependence on \( \mu \) and \(\delta \) more clearly.

Now let us turn to the mean-square recursion of \(\bcw_{e,i}\). First note that \(\rho(\mathcal{J}_{\epsilon}) =
\lambda_2(A) < 1\). Since \(\mathcal{J}_{\epsilon}\) has a Jordan structure, this means that we can chose \( \epsilon \)
small enough, such that \( \|\mathcal{J}_{\epsilon}\|_2=\rho(\mathcal{J}_{\epsilon}^{\T}\mathcal{J}_{\epsilon})\le
\|\mathcal{J}_{\epsilon}^{\T}\mathcal{J}_{\epsilon}\|_{\infty}<1\). Then,
\begin{align}
  \phantom{=}&\: \E \|\bcw_{e,i}\|^2\notag \\
  =&\: \E \Big \|\mathcal{J}_{\epsilon}^{\T} \bcw_{e,i-1} + \mu \mathcal{J}_{\epsilon}^{\T} \mathcal{V}_{R}^{\T}\big( \boldsymbol{t}_{i-1} + \boldsymbol{u}_{i-1} + \s_i  - \E \s_i \notag \\
  &\: \ \ \ \  + \E \s_i - {g}(\mathds{1} \otimes \w_{c,i-1})-{r}(\mathds{1} \otimes \w_{c,i-1})\big)\Big \|^2\notag \\
  \stackrel{(a)}{=}&\: \E \Big \|\mathcal{J}_{\epsilon}^{\T} \bcw_{e,i-1} + \mu \mathcal{J}_{\epsilon}^{\T} \mathcal{V}_{R}^{\T}\big( \boldsymbol{t}_{i-1} + \boldsymbol{u}_{i-1} + \E \s_i \notag \\
  &\: \ \ \ \ -{g}(\mathds{1} \otimes \w_{c,i-1})-{r}(\mathds{1} \otimes \w_{c,i-1})\big)\Big \|^2\notag \\
  &\: + \mu^2\E \left \|\mathcal{J}_{\epsilon}^{\T} \mathcal{V}_{R}^{\T} \left( \s_i - \E \s_i \right)\right \|^2\notag \\
  \stackrel{(b)}{\le}&\: \frac{1}{\|\mathcal{J}_{\epsilon}\|} \E \left \|\mathcal{J}_{\epsilon}^{\T} \bcw_{e,i-1}\right \|^2 + \frac{\mu^2}{1-\|\mathcal{J}_{\epsilon}\|} \notag \\
  &\: \times \E \Big \|\mathcal{J}_{\epsilon}^{\T} \mathcal{V}_{R}^{\T}\big( \boldsymbol{t}_{i-1} + \boldsymbol{u}_{i-1} + \E \s_i \notag \\
  &\: \ \ \ \ \ \ \ \ - {g}(\mathds{1} \otimes \w_{c,i-1})-{r}(\mathds{1} \otimes \w_{c,i-1})\big)\Big \|^2 \notag \\
  &\: + \mu^2\E \left \|\mathcal{J}_{\epsilon}^{\T} \mathcal{V}_{R}^{\T} \left( \s_i - \E \s_i \right)\right \|^2\notag \\
  \stackrel{(c)}{\le}&\: {\|\mathcal{J}_{\epsilon}\|} \E \left \| \bcw_{e,i-1}\right \|^2 + \frac{\mu^2\|\mathcal{J}_{\epsilon}\|^2 \|\mathcal{V}_{R}\|^2}{1-\|\mathcal{J}_{\epsilon}\|} \notag \\
  &\: \times \E \Big \|\boldsymbol{t}_{i-1} + \boldsymbol{u}_{i-1} + \E \s_i \notag \\
  &\: \ \ \ \ \ \ \ \ - {g}(\mathds{1} \otimes \w_{c,i-1})-{r}(\mathds{1} \otimes \w_{c,i-1})\Big \|^2 \notag \\
  &\: + \mu^2\|\mathcal{J}_{\epsilon}\|^2 \|\mathcal{V}_{R}\|^2\E \left \|\s_i - \E \s_i\right \|^2\notag \\
  \stackrel{(d)}{\le}&\: {\|\mathcal{J}_{\epsilon}\|} \E \left \|\bcw_{e,i-1}\right \|^2 + \frac{25 \mu^2\|\mathcal{J}_{\epsilon}\|^2 \|\mathcal{V}_{R}\|^2}{1-\|\mathcal{J}_{\epsilon}\|} \notag \\
  &\: \times \Big(\E \left \|\boldsymbol{t}_{i-1}\right \|^2 + \E \left \|\boldsymbol{u}_{i-1}\right \|^2 + \E \left \|\E \s^p_i\right \|^2 \notag \\
  &\: \ \ \ \ \  + \E \left \|{g}(\mathds{1} \otimes \w_{c,i-1})\right \|^2 +\E \left \|{r}(\mathds{1} \otimes \w_{c,i-1})\right \|^2\Big) \notag \\
  &\: + \mu^2\|\mathcal{J}_{\epsilon}\|^2 \|\mathcal{V}_{R}\|^2\E \left \|\s_i - \E \s_i\right \|^2\notag \\
  \stackrel{(e)}{=}&\: {\|\mathcal{J}_{\epsilon}\|} \E \left \|\bcw_{e,i-1}\right \|^2 + \frac{25 \mu^2\|\mathcal{J}_{\epsilon}\|^2 \|\mathcal{V}_{R}\|^2}{1-\|\mathcal{J}_{\epsilon}\|} \notag \\
  &\: \times \mathds{1}^{\T}\Big(P[\boldsymbol{t}_{i-1}] + P[\boldsymbol{u}_{i-1}] + P[\E \s_i] \notag \\
  &\: \ \ \ \ \ \ \ \ \ + P[{g}(\mathds{1} \otimes \w_{c,i-1})]+P[{r}(\mathds{1} \otimes \w_{c,i-1})]\Big) \notag \\
  &\: + \mu^2\|\mathcal{J}_{\epsilon}\|^2 \|\mathcal{V}_{R}\|^2 \mathds{1}^{\T}P[\s_i - \E \s_i]\notag \\
  \stackrel{(f)}{\le}&\: {\|\mathcal{J}_{\epsilon}\|} \E \left \|\bcw_{e,i-1}\right \|^2 + \frac{25 \mu^2\|\mathcal{J}_{\epsilon}\|^2 \|\mathcal{V}_{R}\|^2}{1-\|\mathcal{J}_{\epsilon}\|} \notag \\
  &\: \times \mathds{1}^{\T}\Big( \Big( 2 \lambda_U^2 + \frac{1+\mu^2}{\delta^2}\Big) \nu^2 \mathds{1}\mathds{1}^{\T} P[\bcw_{e,i-1}]\notag \\
  &\: \ \ \ \ \ \ \ \ \ + \frac{\mu^2}{\delta^2} \big(3\lambda_U^2 \nu^2 \mathds{1}\mathds{1}^{\T}P[\bcw_{e,i-1}] + 3\lambda_U^2 P[\mathds{1}\otimes \widetilde{\w}_{c,i-1}] \notag \\
  &\: \ \ \ \ \ \ \ \ \ \ \ \ \ \ \ \ \ \ + 3 P[g(\mathds{1} \otimes w_{\delta}^o)]\big) \notag \\
  &\: \ \ \ \ \ \ \ \ \ + 3\beta^2 \frac{\mu^2}{\delta^2} P[\mathds{1}\otimes \widetilde{\w}_{c,i-1}] + 3 \beta^2 \frac{\mu^2}{\delta^2} \nu^2 \mathds{1}\mathds{1}^{\T}P[\bcw_{e,i-1}] \notag \\
  &\: \ \ \ \ \ \ \ \ \ + 3\beta^2 \frac{\mu^2}{\delta^2}P[\mathds{1}\otimes w_{\delta}^o] + \frac{\mu^2}{\delta^2}\sigma^2 \mathds{1}\notag \\
  &\: \ \ \ \ \ \ \ \ \ + 2\lambda_U^2 P[\mathds{1}\otimes \widetilde{\w}_{c,i-1}] + 2P[g(\mathds{1}\otimes w_{\delta}^o)] \notag \\
  &\: \ \ \ \ \ \ \ \ \ + \frac{2}{\delta^2} P[\mathds{1}\otimes \widetilde{\w}_{c,i-1}] + 2P[r(\mathds{1}\otimes w_{\delta}^o)]\Big)\notag \\
  &\: +\mu^2\|\mathcal{J}_{\epsilon}\|^2 \|\mathcal{V}_{R}\|^2 \mathds{1}^{\T}\Big(3\beta^2 P[\mathds{1}\otimes \widetilde{\w}_{c,i-1}] \notag \\
  &\: \ \ \ \ \ \ \ \ \ \ \ \ + 3 \beta^2 \nu^2 \mathds{1}\mathds{1}^{\T}P[\bcw_{e,i-1}] + 3\beta^2 P[\mathds{1}\otimes w_{\delta}^o] + \sigma^2 \mathds{1}\Big)\notag \\
  \stackrel{(g)}{=}&\: {\|\mathcal{J}_{\epsilon}\|} \E \left \|\bcw_{e,i-1}\right \|^2 + \frac{25 \mu^2\|\mathcal{J}_{\epsilon}\|^2 \|\mathcal{V}_{R}\|^2}{1-\|\mathcal{J}_{\epsilon}\|} \notag \\
  &\: \times \Big( \Big( 2 \lambda_U^2 + \frac{1 + \mu^2}{\delta^2}\Big) \nu^2 N \E \|\bcw_{e,i-1}\|^2\notag \\
  &\: \ \ \ \ \ + \frac{\mu^2}{\delta^2} \big(3\lambda_U^2 \nu^2 N\E \|\bcw_{e,i-1}\|^2 + 3\lambda_U^2 N \E \|\widetilde{\w}_{c,i-1}\|^2 \notag \\
  &\: \ \ \ \ \ \ \ \ \ \ \ \ \ + 3 \|g(\mathds{1}\otimes w_{\delta}^o)\|^2\big) \notag \\
  &\: \ \ \ \ \ + 3\beta^2 \frac{\mu^2}{\delta^2} N\E \|\widetilde{\w}_{c,i-1}\|^2 + 3 \beta^2 \frac{\mu^2}{\delta^2} \nu^2 N\E \|\bcw_{e,i-1}\|^2 \notag \\
  &\: \ \ \ \ \ + 3\beta^2 \frac{\mu^2}{\delta^2}N \|w_{\delta}^o\|^2 + \frac{\mu^2}{\delta^2}N\sigma^2\notag \\
  &\: \ \ \ \ \ + 2\lambda_U^2 N\E \|\widetilde{\w}_{c,i-1}\|^2 + 2\|g(\mathds{1}\otimes w_{\delta}^o)\|^2 \notag \\
  &\: \ \ \ \ \ + \frac{2}{\delta^2}N \E \|\widetilde{\w}_{c,i-1}\|^2 + 2\|r(\mathds{1}\otimes w_{\delta}^o)\|^2\Big)\notag \\
  &\: + \mu^2\|\mathcal{J}_{\epsilon}\|^2 \|\mathcal{V}_{R}\|^2 \Big(3\beta^2 N\E \|\widetilde{\w}_{c,i-1}\|^2 \notag \\
  &\: \ \ \ \ \ \ \ \ \ \ \ \ \ \ \ \ \ \ \ \ \ \ \ +3 \beta^2 \nu^2 N\E \|\bcw_{e,i-1}\|^2 \notag \\
  &\: \ \ \ \ \ \ \ \ \ \ \ \ \ \ \ \ \ \ \ \ \ \ \ + 3\beta^2 N\|w_{\delta}^o\|^2 + N\sigma^2 \Big)\notag \\
  =&\: \Bigg({\|\mathcal{J}_{\epsilon}\|} + \mu^2\nu^2 N \|\mathcal{J}_{\epsilon}\|^2 \|\mathcal{V}_{R}\|^2 \bigg( \frac{25}{1-\|\mathcal{J}_{\epsilon}\|} \Big(2\lambda_U^2 \notag \\
  &\: \ \ \ \ \ + \frac{1+\mu^2}{\delta^2} + 3 \frac{\mu^2}{\delta^2} ( \lambda_U^2 + \beta^2 )\Big)+3\beta^2\bigg)\Bigg) \E \left \|\bcw_{e,i-1}\right \|^2 \notag \\
  &\: +\mu^2 N \|\mathcal{J}_{\epsilon}\|^2 \|\mathcal{V}_{R}\|^2 \Bigg(\frac{25}{1-\|\mathcal{J}_{\epsilon}\|}\bigg(3 \lambda_U^2 \frac{\mu^2}{\delta^2} + 3\beta^2\frac{\mu^2}{\delta^2} \notag \\
  &\: \ \ \ \ \ + 2\lambda_U^2 + \frac{2}{\delta^2}\bigg) + 3 \beta^2\Bigg)\E \|\widetilde{\w}_{c,i-1}\|^2 \notag \\
  &\:+ \frac{25\mu^2\|\mathcal{J}_{\epsilon}\|^2 \|\mathcal{V}_{R}\|^2}{1-\|\mathcal{J}_{\epsilon}\|}\Bigg(\bigg(2+3 \frac{\mu^2}{\delta^2}\bigg)\|g(\mathds{1}\otimes w_{\delta}^o)\|^2 \notag \\
  &\: \ \ \ \ \ \ \ \ + 3\beta^2 \frac{\mu^2}{\delta^2}N \|w_{\delta}^o\|^2 + \frac{\mu^2}{\delta^2}N\sigma^2+ 2\|r(\mathds{1}\otimes w_{\delta}^o)\|^2\Bigg)\notag \\
  &\:+\mu^2\|\mathcal{J}_{\epsilon}\|^2 \|\mathcal{V}_{R}\|^2 \Big( 3\beta^2 N\|w_{\delta}^o\|^2 + N\sigma^2 \Big)\notag \\
  =&\: \Bigg({\|\mathcal{J}_{\epsilon}\|} + \frac{\mu^2}{\delta^2}\frac{25\nu^2 N \|\mathcal{J}_{\epsilon}\|^2 \|\mathcal{V}_{R}\|^2}{1-\|\mathcal{J}_{\epsilon}\|} \notag \\
  &\: \ \ + \mu^2 \frac{25\nu^2 N \|\mathcal{J}_{\epsilon}\|^2 \|\mathcal{V}_{R}\|^2}{1-\|\mathcal{J}_{\epsilon}\|}\left(2 \lambda_U^2 + \frac{1-\|\mathcal{J}_{\epsilon}\|}{25} 3 \beta^2\right) \notag \\
  &\: \ \ + \frac{\mu^4}{\delta^2} \frac{25\nu^2 N \|\mathcal{J}_{\epsilon}\|^2 \|\mathcal{V}_{R}\|^2}{1-\|\mathcal{J}_{\epsilon}\|} \left( 1+3\beta^2 +3\lambda_U^2 \right) \Bigg) \notag \\
  &\: \ \ \ \ \times \E \left \|\bcw_{e,i-1}\right \|^2 \notag \\
  &\: +\Bigg(\frac{\mu^2}{\delta^2} \frac{50 N \|\mathcal{J}_{\epsilon}\|^2 \|\mathcal{V}_{R}\|^2}{1-\|\mathcal{J}_{\epsilon}\|} \notag \\
  &\: \ \ \ \ \ \ + \mu^2 \frac{25 N \|\mathcal{J}_{\epsilon}\|^2 \|\mathcal{V}_{R}\|^2}{1-\|\mathcal{J}_{\epsilon}\|} \bigg( 2\lambda_U^2 + \frac{1-\|\mathcal{J}_{\epsilon}\|}{25}3 \beta^2\bigg) \notag \\
  &\: \ \ \ \ \ \ + \frac{\mu^4}{\delta^2} \frac{75 N \|\mathcal{J}_{\epsilon}\|^2 \|\mathcal{V}_{R}\|^2}{1-\|\mathcal{J}_{\epsilon}\|} \big( \lambda_U^2 + \beta^2 \big) \Bigg)\E \|\widetilde{\w}_{c,i-1}\|^2 \notag \\
  &\:+ \frac{25\mu^2\|\mathcal{J}_{\epsilon}\|^2 \|\mathcal{V}_{R}\|^2}{1-\|\mathcal{J}_{\epsilon}\|}\Bigg(\bigg(2+3 \frac{\mu^2}{\delta^2}\bigg)\|g(\mathds{1}\otimes w_{\delta}^o)\|^2 \notag \\
  &\: \ \ \ \ \ \ \ \ \ \ + 3\beta^2 \frac{\mu^2}{\delta^2}N \|w_{\delta}^o\|^2 + \frac{\mu^2}{\delta^2}N\sigma^2+ 2\|r(\mathds{1}\otimes w_{\delta}^o)\|^2\Bigg)\notag \\
  &\:+\mu^2\|\mathcal{J}_{\epsilon}\|^2 \|\mathcal{V}_{R}\|^2 \Big( 3\beta^2 N\|w_{\delta}^o\|^2 + N\sigma^2 \Big)
\end{align}
In step (a), cross-terms are eliminated because \(\E \left \{ \s_i - \E \s_i\right \}=0\). Step (b) is due to \( \|
\mathcal{J}_{\epsilon}\| < 1\) and Jensen's inequality, (c) is due to the sub-multiplicative property of norms, (d)
follows from Jensen's inequality, and \((e)\) is due to \(\mathds{1}^{\T} P[\boldsymbol{x}] = \E \|x\|^2\). The bounds
from Lemma~\ref{lem:pertbounds} are used in (f) and (g) is due to \(\mathds{1}^{\T} P[\boldsymbol{x}]=\E
\|\boldsymbol{x}\|^2\) for \(\boldsymbol{x} \in \mathds{R}^{MN}\) and \( \mathds{1}^{\T}
P[\mathds{1}\otimes\boldsymbol{y}]=N\cdot \E \|\boldsymbol{y}\|^2\) for \(\boldsymbol{y} \in \mathds{R}^{M}\).

\section{Proof of Lemma~\ref{lem:lemma8}}\label{ap:lemma8}
\noindent For \(\delta = \mu^{\frac{1}{2}-\kappa}\) and small step-sizes \( \mu \),
\begin{equation}
  \Gamma= \begin{bmatrix} 1 - \mu \lambda_L + O(\mu^2) & O(\mu^{2\kappa})\\
    O(\mu^{1+2\kappa}) & \|\mathcal{J}_{\epsilon}\|+O(\mu^{1+2\kappa})
  \end{bmatrix}\\
\end{equation}
so that
\begin{equation}
  \| \Gamma \|_1 = \max\left \{ 1-\mu \lambda_L + O(\mu^{1+2\kappa}), \| \mathcal{J}_{\epsilon} \| + O(\mu^{2\kappa}) \right \} < 1
\end{equation}
for small enough \(\mu \). Since \( \rho(\Gamma) \le {\| \Gamma \|}_1 < 1 \), \( \Gamma \) is stable. It is also invertible and we obtain
\begin{equation}
  \underset{i \to \infty}{\limsup} \: \begin{bmatrix} \E \|\widetilde{\w}_{c,i}\|^2\\
  \E \|\bcw_{e,i}\|^2 \end{bmatrix} \preceq {\left( I - \Gamma\right)}^{-1} \begin{bmatrix} O(\mu^{2})\\
  O(\mu^{1+2\kappa}) \end{bmatrix}
\end{equation}
Using the matrix inversion lemma, we have
\ifarx
\begin{align}
  {\left( I - \Gamma\right)}^{-1} =& \begin{bmatrix} \mu \lambda_L - O(\mu^2) & -O(\mu^{2\kappa})\\ -O(\mu^{1+2\kappa}) &
    1-\|\mathcal{J}_{\epsilon}\|-O(\mu^{1+2\kappa}) \end{bmatrix}^{-1}\notag \\
    =& \begin{bmatrix} O(\mu) & -O(\mu^{2\kappa})\\ -O(\mu^{1+2\kappa}) & O(1)
  \end{bmatrix}^{-1}\notag \\ =& \begin{bmatrix} O(\mu^{-1}) & O(\mu^{-1 +2\kappa})\\ O(\mu^{2\kappa}) & O(1) \end{bmatrix}
\end{align}
\else
\begin{align}
  {\left( I - \Gamma\right)}^{-1} =& \begin{bmatrix} \mu \lambda_L - O(\mu^2) & -O(\mu^{2\kappa})\\ -O(\mu^{1+2\kappa}) &
    1-\|\mathcal{J}_{\epsilon}\|-O(\mu^{1+2\kappa}) \end{bmatrix}^{-1}\notag \\
    =& \begin{bmatrix} O(\mu^{-1}) & O(\mu^{-1 +2\kappa})\\ O(\mu^{2\kappa}) & O(1) \end{bmatrix}
\end{align}
\fi
The result follows after multiplication and cancellation.

\bibliographystyle{IEEEbib}
\bibliography{prox_transactions}

\begin{IEEEbiography}[{\includegraphics[width=1in,height=1.25in,clip,keepaspectratio]{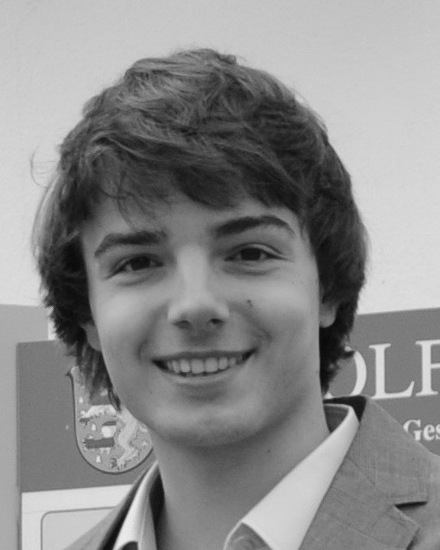}}]{Stefan Vlaski}
  (S'13) received the B.Sc. degree in Electrical Engineering from Technical University Darmstadt, Germany in 2013 and the M.S. in Electrical Engineering as well as Ph.D. in Electrical and Computer Engineering from the University of California, Los Angeles in 2014 and 2019 respectively. He is currently a postdoctoral researcher at the Adaptive Systems Laboratory, EPFL, Switzerland. His research interests are in machine learning, signal processing, and optimization. His current focus is on the development and study of learning algorithms with a particular emphasis on adaptive and decentralized solutions. He was a recipient of the German National Scholarship at TU Darmstadt and the Graduate Division Fellowship at UCLA.
\end{IEEEbiography}
\begin{IEEEbiography}[{\includegraphics[width=1in,height=1.25in,clip,keepaspectratio]{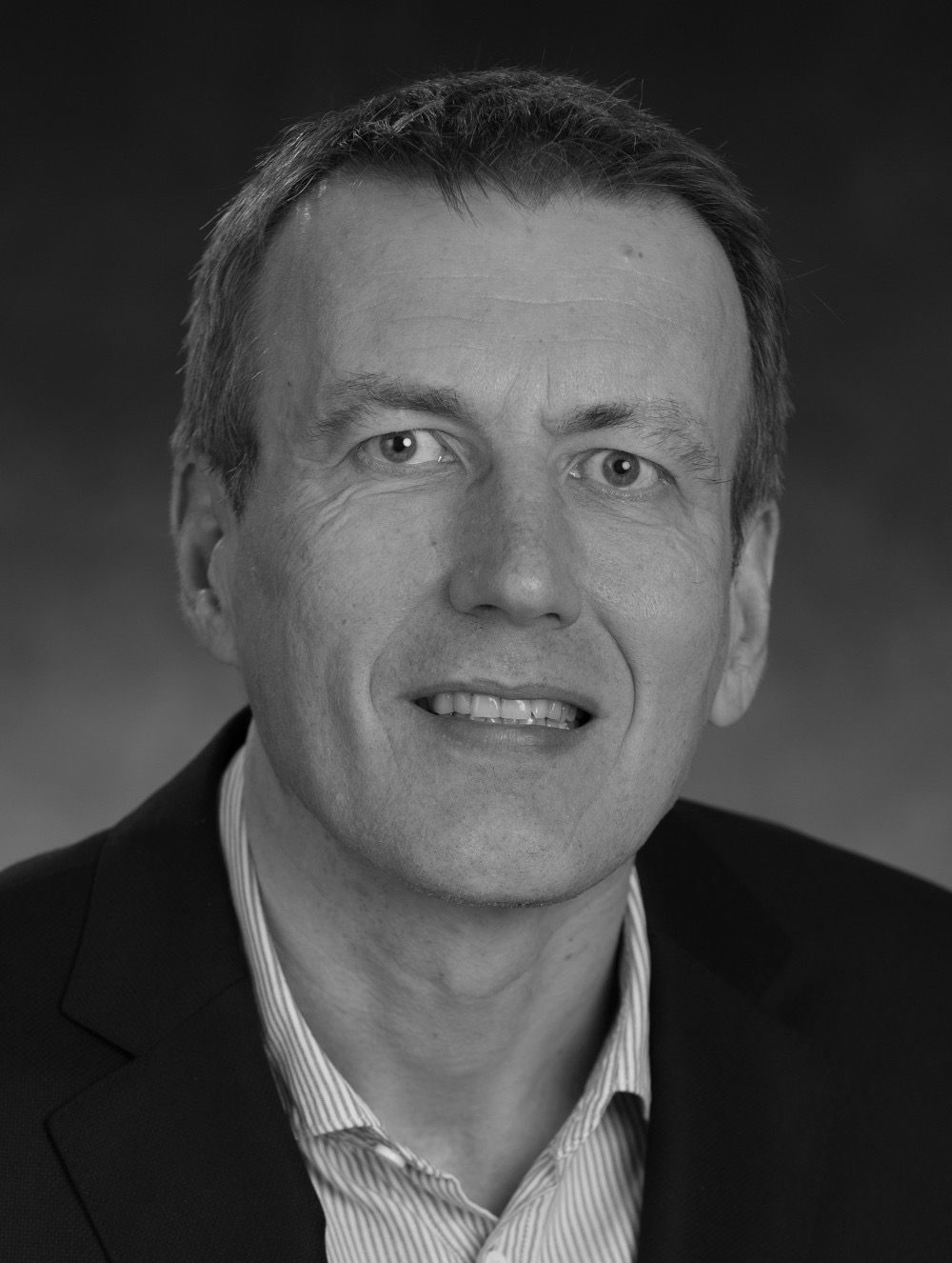}}]{Lieven Vandenberghe}
  is Professor in the Electrical and Computer Engineering Department at UCLA, with a joint appointment in the Department of Mathematics. He obtained the Electrical Engineering degree and a Ph.D. in Applied Sciences from KU Leuven in 1987 and 1992, respectively. He joined UCLA in 1997, following postdoctoral appointments at KU Leuven and Stanford University. Vandenberghe's main research interests are in optimization, systems and control, and signal processing. He is co-author with Stephen Boyd of the books Convex Optimization (2004) and Introduction to Applied Linear Algebra: Vectors, Matrices, and Least Squares (2018). He received the IEEE Guillemin–Cauer Award in 1990, the Robert Stock Award for Exact and Biomedical Sciences at K.U. Leuven in 1993, an NSF CAREER Award in 1998, and the TRW Excellence in Teaching Award of the UCLA Henry Samueli School of Engineering and Applied Science in 2002. He is an Associate Editor of the journals Mathematical Programming and Mathematics of Operations Research.
\end{IEEEbiography}
\begin{IEEEbiography}[{\includegraphics[width=1in,height=1.25in,clip,keepaspectratio]{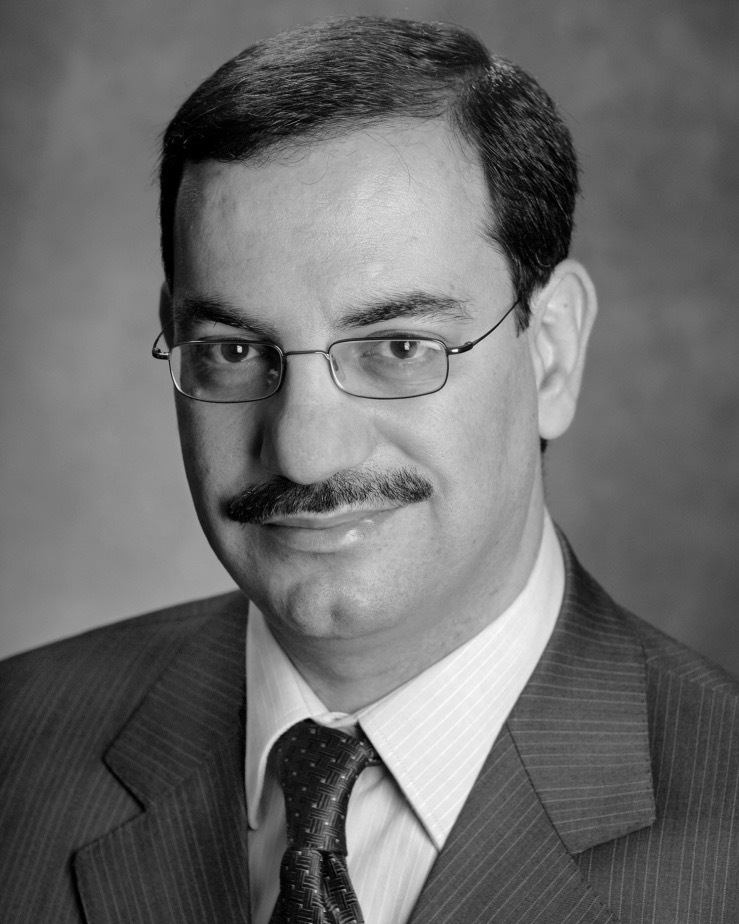}}]{Ali H. Sayed}
  (S'90-M'92-SM'99-F'01) is Dean of Engineering at EPFL, Switzerland, and principal investigator of the Adaptive Systems Laboratory. He has served as distinguished professor and former chairman of electrical engineering at UCLA. An author of over 530 scholarly publications and six books, his research involves several areas including adaptation and learning theories, data and network sciences, statistical inference, multi-agent systems, and optimization. His work has been recognized with several major awards. He is serving as President of the IEEE Signal Processing Society. He is recognized as a Highly Cited Researcher, and is a member of the US National Academy of Engineering.
\end{IEEEbiography}
\vfill

\end{document}